# A SEMIPARAMETRIC MODEL FOR CLUSTER DATA

By Wenyang Zhang, Jianqing Fan[1] and Yan Sun[2]

*University of Bath, Princeton University and Shanghai University of Finance and Economic*

In the analysis of cluster data, the regression coefficients are frequently assumed to be the same across all clusters. This hampers the ability to study the varying impacts of factors on each cluster. In this paper, a semiparametric model is introduced to account for varying impacts of factors over clusters by using cluster-level covariates. It achieves the parsimony of parametrization and allows the explorations of nonlinear interactions. The random effect in the semiparametric model also accounts for within-cluster correlation. Local, linear-based estimation procedure is proposed for estimating functional coefficients, residual variance and within-cluster correlation matrix. The asymptotic properties of the proposed estimators are established, and the method for constructing simultaneous confidence bands are proposed and studied. In addition, relevant hypothesis testing problems are addressed. Simulation studies are carried out to demonstrate the methodological power of the proposed methods in the finite sample. The proposed model and methods are used to analyse the second birth interval in Bangladesh, leading to some interesting findings.

## 1. Introduction.

1.1. *Preamble.* Longitudinal data analysis has attracted considerable attention in the literature. For longitudinal data, the data from the same cluster are dependent with each other. As far as modeling is concerned, this within-cluster dependency is usually accounted by random cluster effects and modeled by a within-cluster correlation matrix. The within-cluster correlation matrix plays a very important role in longitudinal data analysis, as

Received July 2008; revised October 2008.
[1]Supported in part by the NIH Grant R01-GM07261 and NSF Grants DMS-07-04337 and DMS-07-14554.
[2]Supported by National Science Foundation of China (Grant 10801093).
*AMS 2000 subject classifications.* Primary 62G08; secondary 62G10, 62G15.
*Key words and phrases.* Varying-coefficient models, local linear modeling, cluster level variable, cluster effect.

it can be used to improve the efficiency of the estimation. Actually, most of the existing literature is devoted to addressing how to make use of within-cluster correlation matrix to improve the estimation for unknown parameters or functions.

The methodology for parametric based longitudinal data analysis is quite mature (see, e.g., Diggle, Heagerty, Liang and Zeger [6] and the references therein). The situation with nonparametric based longitudinal data analysis is very different. One of the main difficulties is how to incorporate the within cluster correlation structure into the estimation procedure. Lin and Carroll [18] recommend that we ignore the within-cluster correlation when kernel smoothing is employed. Welsh, Lin and Carroll [27] investigate the possibility of using weighted least squares based on the within-cluster correlation structure when spline smoothing is used. They suggest that the weighted least squares estimator based on the true within-cluster correlation structure works better than the estimator based on working independence when spline smoothing is used. Other literature about nonparametric longitudinal regression includes Zeger and Diggle [31], Brumback and Rice [2], Hoover et al. [14], Wu et al. [28], Martinussen and Scheike [21], Chiang et al. [3], Huang et al. [15], Wang [25], Fan and Li [8], Chiou and Müller [4], Wang, Carroll and Lin [26], Qu and Li [23], Lin and Carroll [19], Sun et al. [24] and Fan and Wu [11], among others.

Most of the literature assumes that the regression parameters or functions are the same across all clusters. However, when the regression effects of some particular clusters are of interest, it is unreasonable to assume the regression parameters or functions are the same across all clusters. The interactions of regression effects with clusters are of interest. A naive method to address this is to let the regression coefficients or functions vary freely over clusters. However, this naive method will not be parsimonious, particularly when the number of clusters is large and the issue of estimability arises. In addition, the within-cluster dependency may be addressed by the random effect. This leads us to model these cluster-dependent regression coefficients or functions by using cluster level variables. It addresses, simultaneously, the parsimony and cluster dependency of modeling.

1.2. *A motivating example.* The data that stimulates this project is from Bangladesh and concerns the second birth interval, which is defined as the duration between the first birth and the second birth. The data comes from the Bangladesh Demographic and Health Survey (BDHS) of 1996–1997 (Mitra et al. [22]), which is a cross-sectional, nationally representative survey of ever-married women aged between 10 and 49. Of interest is how some factors that are commonly found to be associated with contraceptive use in Bangladesh, such as education and religion, affect the second birth interval. The data were collected from different districts (clusters) in the six different divisions of Bangladesh. Bangladesh is divided into six administrative divisions: Barisal, Chittagong, Dhaka, Kulna, Rajshahi and Sylhet. The data from the same cluster are correlated with each other, due to cluster-level factors such as cultural norms and access to family planning programs. Of particularly interest is how the factors affect the second birth interval in some particular clusters. For example, how these factors affect the second birth interval in a rural area in Chittagong division.

Some of the factors of interest are defined on the individual level, such as education, and are called individual level variables. Some of the factors are defined on the cluster level, such as type of region of residence, and are called cluster-level variables. We use $y_{ij}$ to denote the length of the second birth interval of the $j$th woman in the $i$th cluster, $X_{ij}$ to denote the vector of the corresponding individual level variables and $Z_i$ to denote the vector of the cluster level variables.

Frequently, the linear model

$$(1.1) \qquad y_{ij} = X_{ij}^T \mathbf{a} + Z_i^T \boldsymbol{\beta} + \varepsilon_{ij}, \qquad j=1,\ldots,n_i, i=1,\ldots,m$$

would be used to fit the data. The within-cluster dependency would be accounted by $\varepsilon_{ij}$, $j=1,\ldots,n_i$ being correlated. The covariance matrix of $\boldsymbol{\varepsilon}_i = (\varepsilon_{i1},\ldots,\varepsilon_{in_i})^T$ can be incorporated into the estimation procedure.

Model (1.1) would be fine if the interest focuses only on the global impact of the factors. However, if the picture for a particular cluster is of interest, (1.1) would not be adequate. Let's take education as an example. It is evident that the impact of education in the cluster where Muslims predominate would be different from the cluster where Hindus predominate. To take the difference of this kind into account, we may relax the assumption imposed on (1.1) and allow the regression coefficients to vary over clusters. This leads to

$$(1.2) \qquad y_{ij} = X_{ij}^T \mathbf{a}_i + Z_i^T \boldsymbol{\beta} + \varepsilon_{ij}, \qquad j=1,\ldots,n_i, i=1,\ldots,m.$$

While it accounts for the varying impact across clusters, (1.2) is not parsimonious. In fact, (1.2) involves $pm+q$ regression coefficients, where $p$ and $q$ are the dimensions of $X_{ij}$ and $Z_i$, respectively. When the number of clusters $m$ is large, there would be too many unknown parameters in model (1.2) for us to get reasonably accurate estimators. In longitudinal data analysis, we often come across large number of clusters. For example, there are 296 clusters in the second birth interval data set that stimulates this paper. If we use (1.2) to fit the data, we would face $296p+q$ unknown coefficients, and would certainly pay a big price on variances of the resulting estimators.

A sensible approach is to model the factor loadings $\mathbf{a}_i$ by using cluster-level variables. A reasonable model is

$$(1.3) \qquad \begin{cases} y_{ij} = X_{ij}^T \mathbf{a}_i + Z_i^T \boldsymbol{\beta} + \varepsilon_{ij}, & j=1,\ldots,n_i,\ i=1,\ldots,m, \\ \mathbf{a}_i = \boldsymbol{\alpha}_0 + \mathbf{A} Z_i + \mathbf{e}_i, & i=1,\ldots,m, \end{cases}$$

where $\mathbf{A} = (\boldsymbol{\alpha}_1, \ldots, \boldsymbol{\alpha}_q)$, and $\mathbf{e}_i$ $(i = 1, \ldots, m)$ are random effects with mean zero. This achieves, simultaneously, the parsimony and within-cluster dependency, and the cluster-dependent factor loadings are allowed. In fact, the number of unknown coefficients in model (1.3) is $p(q+1) + q$, which is usually much smaller than $pm + q$.

A further extension of model (1.3) is to let the factor loadings vary with time as the society and technology evolve with time. By allowing the impacts varying with time, we come up with the model

$$(1.4) \quad \begin{cases} y_{ij} = X_{ij}^T \mathbf{a}_i(U_{ij}) + Z_i^T \boldsymbol{\beta}(U_{ij}) + \varepsilon_{ij}, \\ \mathbf{a}_i(U_{ij}) = \boldsymbol{\alpha}_0(U_{ij}) + \mathbf{A}(U_{ij}) Z_i + \mathbf{e}_i, \end{cases}$$

where $U_{ij}$ is time, $\mathbf{A}(U_{ij}) = (\boldsymbol{\alpha}_1(U_{ij}), \ldots, \boldsymbol{\alpha}_q(U_{ij}))$. Model (1.4) is the model that we are going to address. It is a kind of varying coefficient model (Xia and Li [30], Fan and Zhang [9, 10], Zhang et al. [32] and Li and Liang [17]). To make (1.4) more mathematically clear and general, from now on, $U_{ij}$ is not necessarily to be time, and it can be any continuous covariate. This allows the nonlinear interaction of individual variables $X_{ij}$ and cluster level variable $Z_i$ with $U_{ij}$. We assume that $\varepsilon_{ij}$ is measurement error with mean 0 and variance $\sigma^2$ and independent of $X_{ij}$, $U_{ij}$ and $Z_i$, and that $\{\mathbf{e}_i\}$ are i.i.d. random effects with mean $\mathbf{0}_{p \times 1}$ and covariance matrix $\Sigma$ and independent of all other random variables. We assume that $\{(X_{ij}^T, U_{ij})^T\}$ are i.i.d., and so are $\{Z_i\}$.

In (1.4), $\boldsymbol{\beta}(\cdot)$, $\boldsymbol{\alpha}_k(\cdot)$, $k = 0, \ldots, q$, are unknown functions to be estimated, and so are $\sigma^2$ and $\Sigma$. Although (1.4) is stimulated by the second birth interval data, the modeling concept and estimation methodology, which this paper aims to explore, are equally applicable to other kinds of data, such as the data obtained from medicine and engineering.

The paper is organized as follows. We begin, in Section 2, with a description of the estimation procedure for the proposed model (1.4). In Section 3, we establish the asymptotic properties of the proposed estimators. Hypothesis test associated with model (1.4) is discussed in Section 4. The performance of the method is assessed by a simulation study in Section 5. In Section 6, we use the proposed model and estimation procedure to analyse the data on the second birth intervals in Bangladesh and explore how the impacts of the factors of interest on the length of second birth intervals in some particular clusters change over time.

**2. Estimation procedure.** In this section, we are going to construct the estimation procedure for the proposed model (1.4). We estimate the unknown functional coefficients first, then $\sigma^2$ and $\Sigma$.

2.1. *Estimation of functional coefficients.* By Taylor's expansion, we have

$$\boldsymbol{\alpha}_k(U_{ij}) \approx \boldsymbol{\alpha}_k(u) + \dot{\boldsymbol{\alpha}}_k(u)(U_{ij} - u), \qquad k = 0, \ldots, q,$$

$$\boldsymbol{\beta}(U_{ij}) \approx \boldsymbol{\beta}(u) + \dot{\boldsymbol{\beta}}(u)(U_{ij} - u),$$

when $U_{ij}$ is in a small neighborhood of $u$, which leads to the local least squares estimation procedure

$$L = \sum_{i=1}^{m} \sum_{j=1}^{n_i} \left( y_{ij} - X_{ij}^T \left[ \sum_{k=0}^{q} \{\mathbf{b}_k + \mathbf{c}_k(U_{ij} - u)\} z_{ik} \right] \right. \tag{2.1}$$

$$\left. - Z_i^T \{\mathbf{J} + \mathbf{d}(U_{ij} - u)\} \right)^2 K_h(U_{ij} - u),$$

where $z_{i0} = 1$, $K_h(\cdot) = K(\cdot/h)/h$, $h$ is a bandwidth, and $K(\cdot)$ is a kernel function. We minimize $L$, with respect to $\mathbf{J}$, $\mathbf{d}$, $\mathbf{b}_k$, $\mathbf{c}_k$, $k = 0, \ldots, q$, to get the minimizer $\hat{\mathbf{J}}$, $\hat{\mathbf{d}}$, $\hat{\mathbf{b}}_k$, $\hat{\mathbf{c}}_k$, $k = 0, \ldots, q$. We use $\hat{\mathbf{b}}_k$ to estimate $\boldsymbol{\alpha}_k(u)$ and $\hat{\mathbf{J}}$ to estimate $\boldsymbol{\beta}(u)$. From now on, we denote $\hat{\mathbf{b}}_k$ by $\hat{\boldsymbol{\alpha}}_k(u)$ and $\hat{\mathbf{J}}$ by $\hat{\boldsymbol{\beta}}(u)$.

Let

$$\mathbf{x}_i = (X_{i1}, \ldots, X_{in_i})^T, \qquad \Gamma_i = (\mathbf{x}_i\{(1, Z_i^T) \otimes I_p\}, \mathbf{1}_{n_i} \otimes Z_i^T),$$

$$\Gamma = (\Gamma_1^T, \ldots, \Gamma_m^T)^T, \qquad \mathbf{X} = (\Gamma, \mathcal{U}_1 \Gamma),$$

$$\mathcal{U}_i = \operatorname{diag}((U_{11} - u)^i, \ldots, (U_{1n_1} - u)^i, \ldots, (U_{m1} - u)^i, \ldots, (U_{mn_m} - u)^i),$$

$$W = \operatorname{diag}(K_h(U_{11} - u), \ldots, K_h(U_{1n_1} - u), \ldots,$$

$$K_h(U_{m1} - u), \ldots, K_h(U_{mn_m} - u)),$$

$$Y = (y_{11}, \ldots, y_{1n_1}, \ldots, y_{m1}, \ldots, y_{mn_m}),$$

where $I_p$ is s size $p$ identity matrix and $\mathbf{1}_d$ is a $d$ dimensional vector with each component being 1. By a simple calculation, we have

$$\hat{\boldsymbol{\alpha}}_k(u) = \mathbf{A}_k Y \qquad (k = 0, \ldots, q), \qquad \hat{\boldsymbol{\beta}}(u) = \mathbf{B} Y, \tag{2.2}$$

where

$$\mathbf{A}_k = (e_{(k+1),(q+1)}^T \otimes I_p, \mathbf{0}_{p \times (q+s)})(\mathbf{X}^T W \mathbf{X})^{-1} \mathbf{X}^T W,$$

$$\mathbf{B} = (\mathbf{0}_{q \times ((q+1)p)}, I_q, \mathbf{0}_{q \times s})(\mathbf{X}^T W \mathbf{X})^{-1} \mathbf{X}^T W$$

with $e_{k,p}$ denoting the unit vector of length $p$ with 1 at position $k$, $\mathbf{0}_{p \times q}$ the size $p \times q$ matrix with all entries 0 and $s = (q+1)p + q$.

In practice, some coefficients in (1.4) are constant. Under such a situation, model (1.4) becomes a semivarying coefficient mixed effects model. An interesting question is how to estimate the constant coefficients. Fan and Zhang [9] studied a varying coefficient model with coefficients having different degrees of smoothness. They proposed a two-stage estimation procedure. Based on their idea, we propose a very simple estimation procedure for the unknown constant coefficients. Suppose that the $j$th component $\alpha_{kj}(u)$ of $\boldsymbol{\alpha}_k(u)$ is a constant that is denoted by $C_{kj}$. We first pretend that $\alpha_{kj}(u)$ is a function and get the estimator $\hat{\alpha}_{kj}(U_{il})$ of $\alpha_{kj}(u)$ at $U_{il}$, $l = 1, \ldots, n_i$, $i = 1, \ldots, m$, by the above estimation procedure. Then, take the average of $\hat{\alpha}_{kj}(U_{il})$ over $l = 1, \ldots, n_i$, $i = 1, \ldots, m$. This average

$$(2.3) \qquad \hat{C}_{kj} = \frac{1}{n} \sum_{i=1}^{m} \sum_{l=1}^{n_i} \hat{\alpha}_{kj}(U_{il}), \qquad n = \sum_{i=1}^{m} n_i$$

is our estimator of $C_{kj}$. We will show, later, that the convergence rate of this estimator is of order $O_P(n^{-1/2})$, when the bandwidth is properly selected. This provides a simple method for estimating the constant coefficients.

A more efficient estimate for the constant coefficient can be obtained by using the profile likelihood method (see, e.g., Lam and Fan [16]). For simplicity, we do not pursue this further.

2.2. *Estimation of $\sigma^2$ and $\Sigma$*. Let $\tilde{\mathbf{a}}_i(U_{ij}) = \boldsymbol{\alpha}_0(U_{ij}) + \sum_{k=1}^{q} \boldsymbol{\alpha}_k(U_{ij}) z_{ik}$ and $\mathbf{r}_i = (r_{i1}, \ldots, r_{in_i})^T$ where $r_{ij} = y_{ij} - X_{ij}^T \tilde{\mathbf{a}}_i(U_{ij}) - Z_i^T \boldsymbol{\beta}(U_{ij})$. Correspondingly, let $\hat{\mathbf{a}}_i(\cdot)$ and $\hat{\mathbf{r}}_i$ be their substitution estimators. Set

$$\mathbf{x}_i = (X_{i1}, \ldots, X_{in_i})^T \quad \text{and} \quad P_i = \mathbf{x}_i(\mathbf{x}_i^T \mathbf{x}_i)^{-1} \mathbf{x}_i^T.$$

For each given $i$, we have the linear model

$$(2.4) \qquad \mathbf{r}_i = \mathbf{x}_i \mathbf{e}_i + \boldsymbol{\varepsilon}_i, \qquad \boldsymbol{\varepsilon}_i = (\varepsilon_{i1}, \ldots, \varepsilon_{in_i})^T.$$

The residual sum of squares of this linear model

$$\text{rss}_i = \mathbf{r}_i^T (I_{n_i} - P_i) \mathbf{r}_i$$

would be the raw material for estimating $\sigma^2$. The degree of freedom of $\text{rss}_i$ is $n_i - p$. Let $\text{RSS}_i$ be $\text{rss}_i$ with $\mathbf{r}_i$ replaced by $\hat{\mathbf{r}}_i$. Pooling all $\{\text{RSS}_i\}$ together leads to the estimator of $\sigma^2$ as

$$\hat{\sigma}^2 = (n - mp)^{-1} \sum_{i=1}^{m} \text{RSS}_i, \qquad n = \sum_{i=1}^{m} n_i.$$

Finally, we estimate $\Sigma$. From (2.4), we have the least squares estimator of $\mathbf{e}_i$ as

$$\tilde{\mathbf{e}}_i = (\mathbf{x}_i^T \mathbf{x}_i)^{-1} \mathbf{x}_i^T \mathbf{r}_i = \mathbf{e}_i + (\mathbf{x}_i^T \mathbf{x}_i)^{-1} \mathbf{x}_i^T \boldsymbol{\varepsilon}_i,$$

which leads to

$$\sum_{i=1}^{m} \tilde{\mathbf{e}}_i \tilde{\mathbf{e}}_i^T = \sum_{i=1}^{m} \mathbf{e}_i \mathbf{e}_i^T + \sum_{i=1}^{m} (\mathbf{x}_i^T \mathbf{x}_i)^{-1} \mathbf{x}_i^T \boldsymbol{\varepsilon}_i \boldsymbol{\varepsilon}_i^T \mathbf{x}_i (\mathbf{x}_i^T \mathbf{x}_i)^{-1} + \sum_{i=1}^{m} (\mathbf{x}_i^T \mathbf{x}_i)^{-1} \mathbf{x}_i^T \boldsymbol{\varepsilon}_i \mathbf{e}_i^T$$
$$+ \sum_{i=1}^{m} \mathbf{e}_i \boldsymbol{\varepsilon}_i^T \mathbf{x}_i (\mathbf{x}_i^T \mathbf{x}_i)^{-1}.$$

The last two terms are of order $O_P(m^{1/2})$, so they are negligible. Hence,

$$m^{-1} \sum_{i=1}^{m} \mathbf{e}_i \mathbf{e}_i^T \approx m^{-1} \left\{ \sum_{i=1}^{m} \tilde{\mathbf{e}}_i \tilde{\mathbf{e}}_i^T - \sum_{i=1}^{m} (\mathbf{x}_i^T \mathbf{x}_i)^{-1} \mathbf{x}_i^T \boldsymbol{\varepsilon}_i \boldsymbol{\varepsilon}_i^T \mathbf{x}_i (\mathbf{x}_i^T \mathbf{x}_i)^{-1} \right\}$$
$$\approx m^{-1} \left\{ \sum_{i=1}^{m} \tilde{\mathbf{e}}_i \tilde{\mathbf{e}}_i^T - \sigma^2 \sum_{i=1}^{m} (\mathbf{x}_i^T \mathbf{x}_i)^{-1} \right\}.$$

Therefore, we have

$$(2.5) \qquad \hat{\Sigma} = m^{-1} \sum_{i=1}^{m} \hat{\mathbf{e}}_i \hat{\mathbf{e}}_i^T - m^{-1} \hat{\sigma}^2 \sum_{i=1}^{m} (\mathbf{x}_i^T \mathbf{x}_i)^{-1}$$

to estimate $\Sigma$. In (2.5), $\hat{\mathbf{e}}_i$ is $\tilde{\mathbf{e}}_i$ with $\mathbf{r}_i$ replaced by $\hat{\mathbf{r}}_i$.

**3. Asymptotic properties.** In this section, we are going to present the asymptotic properties of the proposed estimators. For any $p \times p$ symmetric matrix $A$, we use $\text{vech}(A)$ to denote the vector consisting of all elements on and below the diagonal of the matrix $A$, and we use $\text{vec}(A)$ to denote the vector by simply stacking the column vectors of matrix $A$ below one another. Obviously, there exists a unique $p^2 \times p(p+1)/2$ matrix $R_p$ such that $\text{vec}(A) = R_p \text{vech}(A)$.

We first introduce some notation. For any function or function vector $g(\cdot)$, we use $\dot{g}(\cdot)$ and $\ddot{g}(\cdot)$ to denote its first and second derivatives, respectively. We use $\mathcal{D}$ to denote the collections of all individual and cluster level covariates. Let $\mu_i = \int t^i K(t)\,dt$, $\nu_i = \int t^i K^2(t)\,dt$, and let

$$\Omega(u) = E\{(X^T, Z^T \otimes X^T, Z^T)^T (X^T, Z^T \otimes X^T, Z^T) | U = u\}$$
$$= \begin{pmatrix} \Omega_1(u) & \Omega_2(u) \\ \Omega_2(u)^T & \Omega_3(u) \end{pmatrix},$$

where $\Omega_1(u)$ and $\Omega_3(\cdot)$ are, respectively, $p(q+1) \times p(q+1)$ and $q \times q$ submatrix of $\Omega(u)$. Without loss of generality, we will assume that $\mu_0 = 1$.

Our main asymptotic results are presented through the following 6 theorems. We leave the proofs of these theorems to the Appendix. The first two theorems give the asymptotic normality of the estimated functional coefficients.

THEOREM 1. *Under the conditions* (1)–(6) *in the Appendix, when* $nh^5$ *is bounded, for* $k = 0, \ldots, q$, *we have*

$$\sqrt{nhf(u)}\left\{\hat{\boldsymbol{\alpha}}_k(u) - \boldsymbol{\alpha}_k(u) - h^2\frac{\mu_2}{2}\ddot{\boldsymbol{\alpha}}_k(u)\right\}$$

$$\stackrel{D}{\longrightarrow} N_p(\mathbf{0}_{p\times 1},$$

$$\nu_0\{e_{(k+1),(q+1)}^T \otimes I_p\}[\sigma^2\Lambda_1(u)^{-1} + \Theta_1(u)]\{e_{(k+1),(q+1)} \otimes I_p\}),$$

*where*

$$\Lambda_1(u) = \Omega_1(u) - \Omega_2(u)\Omega_3(u)^{-1}\Omega_2(u)^T, \qquad \Theta_1(u) = \Upsilon_1(u)\Xi_1(u)\Upsilon_1(u)^T,$$
$$\Upsilon_1(u) = (\Lambda_1(u)^{-1}, \Lambda_2(u)), \qquad \Lambda_2(u) = -\Lambda_1(u)^{-1}\Omega_2(u)\Omega_3(u)^{-1}$$

*and*

$$\Xi_1(u) = E\{X^T\Sigma X(X^T, Z^T \otimes X^T, Z^T)^T(X^T, Z^T \otimes X^T, Z^T)|U = u\}.$$

THEOREM 2. *Under the conditions* (1)–(6) *in Appendix, when* $nh^5$ *is bounded, we have*

$$\sqrt{nhf(u)}\left\{\hat{\boldsymbol{\beta}}(u) - \boldsymbol{\beta}(u) - h^2\frac{\mu_2}{2}\ddot{\boldsymbol{\beta}}(u)\right\} \stackrel{D}{\longrightarrow} N_q(\mathbf{0}_{q\times 1}, \nu_0[\sigma^2\Lambda_3(u) + \Theta_2(u)]),$$

*where* $\Lambda_3(u) = \Omega_3(u)^{-1} + \Omega_3(u)^{-1}\Omega_2(u)^T\Lambda_1(u)^{-1}\Omega_2(u)\Omega_3(u)^{-1}$, $\Theta_2(u) = \Upsilon_2(u)\Xi_1(u)\Upsilon_2(u)^T$ *and* $\Upsilon_2(u) = (\Lambda_2(u)^T, \Lambda_3(u))$.

We now present the asymptotic normality for the parametric component. To present the asymptotic property of $\sigma^2$, we assume that the following limits exist and are finite: $c_1 = \lim_{m\to\infty} n/(n - mp), c_2 = \lim_{m\to\infty} m/(n - mp)$ and

$$\gamma = \plim_{n\to\infty}(n - mp)^{-1} \sum_{i=1}^{m}\sum_{j=1}^{n_i}[X_{ij}^T(\mathbf{x}_i^T\mathbf{x}_i)^{-1}X_{ij}]^2,$$

where "plim" denotes convergence in probability.

THEOREM 3. *Under the conditions* (1)–(7) *in the Appendix, when* $nh^8 \to 0$, *we have*

$$\sqrt{n}\{\hat{\sigma}^2 - \sigma^2\} \stackrel{D}{\longrightarrow} N(0, 2\sigma^4 c_1(c_1 - 1 - \gamma) + \text{var}(\varepsilon_{11}^2)c_1(2 - c_1 + \gamma)).$$

Theorem 3 suggests the estimator $\hat{\sigma}^2$ is of convergence rate $O_P(n^{-1/2})$, which is the optimal convergence rate of the parametric estimator.

Additional notation is needed for presenting the asymptotic normality of $\Sigma$. Write

$$\Delta_1 = \plim_{m\to\infty} m^{-1} \sum_{i=1}^{m} [(\mathbf{x}_i^T \mathbf{x}_i)^{-1}],$$

$$\Delta_2 = \plim_{m\to\infty} m^{-1} \sum_{i=1}^{m} [(\mathbf{x}_i^T \mathbf{x}_i)^{-1} \otimes (\mathbf{x}_i^T \mathbf{x}_i)^{-1}],$$

$$\Delta_3 = \begin{pmatrix} \Sigma \otimes \Delta_{1(1)} + \Delta_1 \otimes \Sigma_{(1)} \\ \vdots \\ \Sigma \otimes \Delta_{1(p)} + \Delta_1 \otimes \Sigma_{(p)} \end{pmatrix},$$

where $\Delta_{1(r)}, \Sigma_{(r)}$ $(r = 1, \ldots, p)$ denote the $r$th row of $\Delta_1, \Sigma$, respectively. Let

$$\Delta_4 = \plim_{m\to\infty} m^{-1} \sum_{i=1}^{m} \sum_{j=1}^{n_i} [\text{vec}((\mathbf{x}_i^T \mathbf{x}_i)^{-1} X_{ij} X_{ij}^T (\mathbf{x}_i^T \mathbf{x}_i)^{-1})$$
$$\times \text{vec}((\mathbf{x}_i^T \mathbf{x}_i)^{-1} X_{ij} X_{ij}^T (\mathbf{x}_i^T \mathbf{x}_i)^{-1})^T],$$

$$\Delta_5 = \plim_{m\to\infty} m^{-1} \sum_{i=1}^{m} [(\mathbf{x}_i^T \mathbf{x}_i)^{-1} \mathbf{x}_i^T \tilde{P}_i \mathbf{x}_i (\mathbf{x}_i^T \mathbf{x}_i)^{-1}],$$

where $\tilde{P}_i$ is a diagonal matrix generated from the diagonal elements of $P_i$.

THEOREM 4. *Under the conditions* (1)–(7) *in Appendix, when $nh^8 \to 0$, we have*

$$\sqrt{n} \, \text{vech}(\hat{\Sigma} - \Sigma) \xrightarrow{D} N_r(\mathbf{0}_{r\times 1}, (1/c_2 + p)(R_p^T R_p)^{-1} R_p^T \Delta R_p (R_p^T R_p)^{-1}),$$

*where $r = p(p+1)/2$ and*

$$\Delta = E\{(\mathbf{e}_1 \mathbf{e}_1^T) \otimes (\mathbf{e}_1 \mathbf{e}_1^T)\} - \text{vec}(\Sigma) \text{vec}(\Sigma)^T + 2\sigma^4 \Delta_2$$
$$+ \sigma^2 \{\Sigma \otimes \Delta_1 + \Delta_1 \otimes \Sigma + \Delta_3\}$$
$$+ [\text{var}(\varepsilon_{11}^2) - 2\sigma^4]\{\Delta_4 - c_2[2\text{vec}(\Delta_1)\text{vec}(\Delta_1)^T$$
$$- \text{vec}(\Delta_1)\text{vec}(\Delta_2)^T - \text{vec}(\Delta_2)\text{vec}(\Delta_1)^T]\}$$
$$+ \{2\sigma^4 c_1(c_1 - 1 - \gamma) + \text{var}(\varepsilon_{11}^2) c_1(2 - c_1 + \gamma)\} \text{vec}(\Delta_1) \text{vec}(\Delta_1)^T.$$

Theorem 4 shows that the estimator $\hat{\Sigma}$ also achieves the convergence rate of $O_P(n^{-1/2})$.

THEOREM 5. *When $\alpha_{kj}(u)$ is a constant, under the conditions* (1)–(6) *and* (8) *in the Appendix and $nh^4 \to 0$, we have*

$$\sqrt{n}\{\hat{C}_{kj} - C_{kj}\}$$
$$\xrightarrow{D} N(0, e_{kp+j,((q+1)p)}^T \{\sigma^2 E[\Lambda_1(U)^{-1}] + E[\Theta_1(U)] + \Xi_2\} e_{kp+j,((q+1)p)}),$$

*where*

$$\Xi_2 = \plim_{n\to\infty} \frac{1}{n} \sum_{i=1}^{m} \sum_{l=1}^{n_i} \sum_{r=1, r\neq l}^{n_i} \{\Upsilon_1(U_{il}) E[X_{il}^T \Sigma X_{ir} \Gamma_{il} \Gamma_{ir}^T | U_{il}, U_{ir}] \Upsilon_1(U_{ir})^T\},$$

*with* $\Gamma_{ij}^T = (X_{ij}^T, Z_i^T \otimes X_{ij}^T, Z_i^T)$ *being the $j$th row of matrix $\Gamma_i$.*

When the $j$th $(j = 1, \ldots, q)$ component $\beta_j(u)$ of $\boldsymbol{\beta}(u)$ is a constant, Theorem 5 applies to its estimator as well, except that the variance is

$$e_{j,q}^T \{\sigma^2 E[\Lambda_3(U)^{-1}] + E[\Theta_2(U)] + \Xi_3\} e_{j,q},$$

where

$$\Xi_3 = \plim_{n\to\infty} \frac{1}{n} \sum_{i=1}^{m} \sum_{l=1}^{n_i} \sum_{r=1, r\neq l}^{n_i} \{\Upsilon_2(U_{il}) E[X_{il}^T \Sigma X_{ir} \Gamma_{il} \Gamma_{ir}^T | U_{il}, U_{ir}] \Upsilon_2(U_{ir})^T\}.$$

THEOREM 6. *Under the conditions* (1)–(6), (8) *and* (9) *in the Appendix, with $K(t)$ having a compact support $[-c_0, c_0]$ and $h = n^{-\rho}, 1/5 \leq \rho < 1/3$ for all $u \in [a, b]$, we have*

$$P\bigg((-2\log\{h/(b-a)\})^{1/2}$$
$$\times \bigg\{\sup_{u\in[a,b]} \bigg|\frac{\hat{\alpha}_{kj}(u) - \alpha_{kj}(u) - \mathrm{bias}(\hat{\alpha}_{kj}(u)|\mathcal{D})}{[\mathrm{var}\{\hat{\alpha}_{kj}(u)|\mathcal{D}\}]^{1/2}}\bigg| - \omega_n\bigg\} < x\bigg)$$
$$\to \exp\{-2\exp\{-x\}\},$$

*where*

$$\omega_n = (-2\log\{h/(b-a)\})^{1/2}$$
$$+ \frac{1}{(-2\log\{h/(b-a)\})^{1/2}} \bigg[\log \frac{K^2(c_0)}{\nu_0 \pi^{1/2}} + \frac{1}{2}\log\log\{(b-a)/h\}\bigg]$$

*if $K(c_0) \neq 0$ and*

$$\omega_n = (-2\log\{h/(b-a)\})^{1/2} + \frac{1}{(-2\log\{h/(b-a)\})^{1/2}} \log\bigg\{\frac{1}{4\nu_0 \pi} \int (\dot{K}(t))^2 \, dt\bigg\}$$

*if $K(c_0) = 0$, $K(t)$ is absolutely continuous and $K^2(t), (\dot{K}(t))^2$ are integrable on $(-\infty, +\infty)$.*

REMARK. Theorem 6 gives the distribution of the maximum discrepancy between the estimated functional coefficient and the true coefficient. It is the basis for constructing the hypothesis test or confidence band. Theorem 6 also applies to the estimator of any component of $\boldsymbol{\beta}(\cdot)$.

**4. Confidence bands and hypothesis test.** In this section, we will investigate how to construct the confidence bands for the functional coefficients in model (1.4).

For model (1.4), we often wish to know if an estimated functional coefficient is significantly different from zero or if the estimated functional coefficient is really varying. More generally, we wish to test

$$(4.1) \qquad H_0 : \alpha_{kj}(u) = \alpha_0(u) \quad \longleftrightarrow \quad H_1 : \alpha_{kj}(u) \neq \alpha_0(u),$$

where $\alpha_0(u)$ is a specific function. This kind of nonparametric testing problem can be conveniently handled by using the generalized likelihood ratio method (Fan, Zhang and Zhang [12]). Instead, our test statistics will be based on the constructed confidence bands.

As the proposed confidence bands and test statistics involve the estimation of the biases and variances of the proposed estimators of the functional coefficients, we first construct the estimation procedure for the biases and variances. Throughout this paper, for any functional vector $\mathbf{g}(u)$, we use $\mathbf{g}^{(i)}(u)$ to denote the $i$th derivative of $\mathbf{g}(u)$.

4.1. *Estimation for bias and variance.* Following Fan and Gijbels [7], by (2.2), we have, for $k = 0, \ldots, q$, that

$$(4.2) \quad \begin{cases} \operatorname{bias}(\hat{\boldsymbol{\alpha}}_k(u)|\mathcal{D}) = (e^T_{(k+1),(q+1)} \otimes I_p, \mathbf{0}_{p \times (q+s)})(\mathbf{X}^T W \mathbf{X})^{-1} \mathbf{X}^T W \mathbf{R}, \\ \operatorname{bias}(\hat{\boldsymbol{\beta}}(u)|\mathcal{D}) = (\mathbf{0}_{q \times ((q+1)p)}, I_q, \mathbf{0}_{q \times s})(\mathbf{X}^T W \mathbf{X})^{-1} \mathbf{X}^T W \mathbf{R}, \end{cases}$$

where $\mathbf{R} = (R_{11}, \ldots, R_{1n_1}, \ldots, R_{m1}, \ldots, R_{mn_m})^T$ with

$$R_{ij} = X_{ij}^T A_i(U_{ij}) + Z_i^T \{\boldsymbol{\beta}(U_{ij}) - \boldsymbol{\beta}(u) - \dot{\boldsymbol{\beta}}(u)(U_{ij} - u)\},$$

$$A_i(U_{ij}) = \sum_{k=0}^{q} \{\boldsymbol{\alpha}_k(U_{ij}) - \boldsymbol{\alpha}_k(u) - \dot{\boldsymbol{\alpha}}_k(u)(U_{ij} - u)\} z_{ik}, \qquad z_{i0} = 1.$$

By Taylor's expansion, we have

$$\boldsymbol{\alpha}_k(U_{ij}) - \boldsymbol{\alpha}_k(u) - \dot{\boldsymbol{\alpha}}_k(u)(U_{ij} - u)$$
$$\approx 2^{-1} \ddot{\boldsymbol{\alpha}}_k(u)(U_{ij} - u)^2 + 6^{-1} \boldsymbol{\alpha}_k^{(3)}(u)(U_{ij} - u)^3,$$
$$\boldsymbol{\beta}(U_{ij}) - \boldsymbol{\beta}(u) - \dot{\boldsymbol{\beta}}(u)(U_{ij} - u)$$
$$\approx 2^{-1} \ddot{\boldsymbol{\beta}}(u)(U_{ij} - u)^2 + 6^{-1} \boldsymbol{\beta}^{(3)}(u)(U_{ij} - u)^3.$$

This leads to

$$(4.3) \qquad \mathbf{R} \approx (2^{-1} \mathcal{U}_2 \Gamma, 6^{-1} \mathcal{U}_3 \Gamma)(\ddot{\boldsymbol{\theta}}(u)^T, \boldsymbol{\theta}^{(3)}(u)^T)^T,$$

where $\boldsymbol{\theta}(u) = (\boldsymbol{\alpha}_0(u)^T, \ldots, \boldsymbol{\alpha}_p(u)^T, \boldsymbol{\beta}(u)^T)^T$. The estimators of $\ddot{\boldsymbol{\theta}}(u)$ and $\boldsymbol{\theta}^{(3)}(u)$ can be easily obtained by using local cubic fit with an appropriate pilot bandwidth $h^*[= O(n^{-1/7})]$ in Section 2.1, which is optimal for estimating $\ddot{\boldsymbol{\theta}}(u)$. We denote the estimators by $\hat{\ddot{\boldsymbol{\theta}}}(u)$ and $\hat{\boldsymbol{\theta}}^{(3)}(u)$. Substituting them into (4.3), we obtain $\hat{\mathbf{R}}$ and, hence, estimated biases $\widehat{\text{bias}}(\hat{\boldsymbol{\alpha}}_k(u)|\mathcal{D})$ and $\widehat{\text{bias}}(\hat{\boldsymbol{\beta}}(u)|\mathcal{D})$.

We now derive an estimator of variance of $\hat{\boldsymbol{\alpha}}_k(u)$ and $\hat{\boldsymbol{\beta}}(u)$. We notice that the estimators are linear in $Y$, according to (2.2). Hence, we need only to estimate $\text{var}(Y)$. A natural estimator is

$$\widehat{\text{var}}(Y|\mathcal{D}) = \hat{\sigma}^2 I_n + \text{diag}(\mathbf{x}_1 \hat{\Sigma} \mathbf{x}_1^T, \ldots, \mathbf{x}_m \hat{\Sigma} \mathbf{x}_m^T).$$

This, together with (2.2), give us the estimators

$$\widehat{\text{var}}(\hat{\boldsymbol{\alpha}}_k(u)|\mathcal{D}) = \mathbf{A}_k \widehat{\text{var}}(Y|\mathcal{D}) \mathbf{A}_k \quad \text{and} \quad \widehat{\text{var}}(\hat{\boldsymbol{\beta}}(u)|\mathcal{D}) = \mathbf{B}\widehat{\text{var}}(Y|\mathcal{D})\mathbf{B},$$

with $\mathbf{A}_k$ and $\mathbf{B}$ defined in (2.2).

4.2. *Confidence bands and hypothesis testing.* We first state the theorem, based on which the confidence bands and hypothesis tests are constructed.

THEOREM 7. *Under the conditions* (1)–(9) *in the Appendix, with* $K(t)$ *having a compact support* $[-c_0, c_0]$ *and* $h = n^{-\rho}, 1/5 \leq \rho < 1/3$. *Furthermore, if* $\alpha_{kj}^{(3)}(\cdot)$ *is continuous on* $[a, b]$ *and the pilot bandwidth* $h^*$ *is of order* $n^{-1/7}$, *then, for all* $u \in [a, b]$, *we have*

$$P\bigg((-2\log\{h/(b-a)\})^{1/2}$$
$$\times \bigg\{\sup_{u\in[a,b]} \bigg|\frac{\hat{\alpha}_{kj}(u) - \alpha_{kj}(u) - \widehat{\text{bias}}(\hat{\alpha}_{kj}(u)|\mathcal{D})}{[\widehat{\text{var}}\{\hat{\alpha}_{kj}(u)|\mathcal{D}\}]^{1/2}}\bigg| - \omega_n\bigg\} < x\bigg)$$
$$\to \exp\{-2\exp\{-x\}\},$$

*where* $\omega_n$ *is exactly the same as that in Theorem* 6.

REMARK. If the component of $\boldsymbol{\beta}^{(3)}(\cdot)$ is continuous on $[a, b]$, then Theorem 7 applies to its estimator as well.

Based on Theorem 7, the $1-\alpha$ confidence bands of $\alpha_{kj}(u)$ can be easily constructed as

$$(\hat{\alpha}_{kj}(u) - \widehat{\text{bias}}(\hat{\alpha}_{kj}(u)|\mathcal{D}) \pm \Delta_{1,\alpha}(u)),$$

where

$$\Delta_{1,\alpha}(u) = (\omega_n + [\log 2 - \log\{-\log(1-\alpha)\}](-2\log\{h/(b-a)\})^{-1/2})$$
$$\times \{\widehat{\text{var}}(\hat{\alpha}_{kj}(u)|\mathcal{D})\}^{1/2}.$$

It is worthwhile to mention that Xia [29] investigated bias-corrected confidence bands for univariate nonparametric regression.

By Theorem 7, hypothesis (4.1) can be tested by using the test statistic

$$(-2\log\{h/(b-a)\})^{1/2}\bigg\{\sup_{u\in[a,b]}\bigg|\frac{\hat{\alpha}_{kj}(u)-\alpha_0(u)-\widehat{\text{bias}}(\hat{\alpha}_{kj}(u)|\mathcal{D})}{[\widehat{\text{var}}\{\hat{\alpha}_{kj}(u)|\mathcal{D}\}]^{1/2}}\bigg|-\omega_n\bigg\}$$

and rejecting $H_0$ when the test statistic exceeds the asymptotic critical value $c_\alpha = -\log\{-0.5\log(1-\alpha)\}$. Similarly, we may want to ask whether a specific functional coefficient is really varying. This amounts to testing the composite null hypothesis

(4.4) $$H_0: \alpha_{kj}(u) = C_{kj} \quad \longleftrightarrow \quad H_1: \alpha_{kj}(u) \neq C_{kj}.$$

Based on Theorems 5 and 7, we test the problem (4.4) by computing the statistic

$$(-2\log\{h/(b-a)\})^{1/2}\bigg\{\sup_{u\in[a,b]}\bigg|\frac{\hat{\alpha}_{kj}(u)-\hat{C}_{kj}-\widehat{\text{bias}}(\hat{\alpha}_{kj}(u)|\mathcal{D})}{[\widehat{\text{var}}\{\hat{\alpha}_{kj}(u)|\mathcal{D}\}]^{1/2}}\bigg|-\omega_n\bigg\}$$

and rejecting $H_0$ for large values of the test statistic.

**5. Simulation study.** In this section, we are going to use a simulated example to demonstrate how well the proposed estimation method works and examine how much loss one could incur by ignoring the structure of $\mathbf{a}_i(\cdot)$ when the structure holds. We first demonstrate the accuracy of the proposed estimators.

In model (1.4), we take $p = 3$, $q = 2$ and $m = 100$. The cluster sizes $n_i$ ($i = 1, \ldots, m$) are generated by the integer part of $|2\xi| + 6$, $\xi \sim N(0,1)$. The covariants $\{X_{ij}\}$ are independently generated from $N(0, I_p)$, $\{Z_i\}$ are independently generated from $N(0, I_q)$ and $\{U_{ij}\}$ are independently generated from $U(0,1)$. We also set the random effect $\mathbf{e}_i$ following the normal distribution $N(\mathbf{0}_p, \Sigma)$, measurement error $\varepsilon_{ij}$ following normal distribution $N(0, \sigma^2)$ where $\Sigma = (\sigma_{ij}) = 0.5^2 I_p$ and $\sigma = 0.5$. We set $\beta_0(u) = \sin(2\pi u)$ (intercept term), $\boldsymbol{\alpha}_0(u) = (\alpha_{01}(\cdot), \alpha_{02}(\cdot), \alpha_{03}(\cdot))^T$ a vector with each component being $\sin(2\pi u)$, $\boldsymbol{\alpha}_1(u) = (\alpha_{11}(\cdot), \alpha_{12}(\cdot), \alpha_{13}(\cdot))^T$ a vector with each component being $\cos(2\pi u)$, $\boldsymbol{\alpha}_2(u) = (\alpha_{21}(\cdot), \alpha_{22}(\cdot), \alpha_{23}(\cdot))^T$ a vector with each component being $\sin(\pi u)$, $\boldsymbol{\beta}(u) = (\beta_1(\cdot), \beta_2(\cdot))^T$ a vector with each component being $\sin(2\pi u)$.

For any function or functional vector $\mathbf{g}(\cdot)$, if $\hat{\mathbf{g}}(\cdot)$ is an estimator of $\mathbf{g}(\cdot)$, we define the mean integrated squared error (MISE) of $\hat{\mathbf{g}}(\cdot)$ as

$$E\bigg[\int \|\hat{\mathbf{g}}(u) - \mathbf{g}(u)\|^2 \, du\bigg],$$

where $\|\mathbf{b}\|^2 = \mathbf{b}^T \mathbf{b}$ and use it to assess the accuracy of the estimators.

The proposed estimation method is employed to estimate the functional coefficients. The kernel function is taken to be Epanechnikov kernel $0.75(1-t^2)_+$ and bandwidth is taken to be 0.15. The MISEs of the proposed estimators of unknown functions are presented in Table 1, based on 100 simulations. The MSEs of the estimators of $\Sigma$ and $\sigma^2$ are presented in Table 2. From Tables 1 and 2, we can see the proposed estimators are quite accurate.

To visualize the accuracy of the proposed estimators, among the 100 simulations, we single out the one with median performance and plot the estimated functions together with their 95% confidence bands in Figure 1. It shows, again, that the proposed estimators are very accurate.

Now, we turn to examine how much loss one could incur when ignoring the structure of $\mathbf{a}_i(\cdot)$, namely, examining the gain of our model (1.4). We now assume the cluster effects $\mathbf{e}_i = 0$. Given the sizes of the clusters in the above simulated example, many clusters would be too small for one to estimate the corresponding $\mathbf{a}_i(\cdot)$ when ignoring the structure of $\mathbf{a}_i(\cdot)$. So, we now assume all clusters share the same size of 50, and keep the number of clusters 100. Except the change on the sizes of the clusters and the assumption of $\mathbf{e}_i = 0$, all remain the same as the above simulated example.

We use $\text{MISE}_{1,i}$ to denote the MISE of the estimator of $\mathbf{a}_i(\cdot)$ obtained by the proposed estimation method and $\text{MISE}_{2,i}$ to denote the MISE of the estimator obtained without using the structure of $\mathbf{a}_i(\cdot)$; that is, regarding $\mathbf{a}_i(\cdot)$ as a free unknown function and estimating it based on the first part of model (1.4). The ratio

$$\text{RMISE} = \sum_{i=1}^m \text{MISE}_{1,i} \bigg/ \sum_{i=1}^m \text{MISE}_{2,i}$$

is used to assess the loss incurred on the estimation for $\mathbf{a}(\cdot) = (\mathbf{a}_1(\cdot), \ldots, \mathbf{a}_m(\cdot))$ due to ignoring the structure of $\mathbf{a}_i(\cdot)$, $i = 1, \ldots, m$. We compute the RMISEs for different bandwidths and plot the obtained RMISEs against the bandwidths in Figure 2. It is clear that the RMISE is almost 0 when the bandwidth is less than 0.3, and it never goes beyond 0.25, which suggests the loss is significant. Similarly, we define the RMISE for $\boldsymbol{\beta}(\cdot)$ and plot the RMISEs against different bandwidths. It again shows the loss incurred on the estimation for $\boldsymbol{\beta}(\cdot)$ due to ignoring the structure of $\mathbf{a}_i(\cdot)$, and $i=1,\ldots,m$, is still significant.

TABLE 1
*The MISEs of the estimators*

| Estimator | $\hat{\beta}_0(\cdot)$ | $\hat{\alpha}_{01}(\cdot)$ | $\hat{\alpha}_{02}(\cdot)$ | $\hat{\alpha}_{03}(\cdot)$ | $\hat{\alpha}_{11}(\cdot)$ | $\hat{\alpha}_{12}(\cdot)$ |
|---|---|---|---|---|---|---|
| MISE | 0.013 | 0.019 | 0.019 | 0.020 | 0.019 | 0.020 |
| Estimator | $\hat{\alpha}_{13}(\cdot)$ | $\hat{\alpha}_{21}(\cdot)$ | $\hat{\alpha}_{22}(\cdot)$ | $\hat{\alpha}_{23}(\cdot)$ | $\hat{\beta}_1(\cdot)$ | $\hat{\beta}_2(\cdot)$ |
| MISE | 0.019 | 0.016 | 0.017 | 0.016 | 0.014 | 0.014 |

TABLE 2
*The MSEs of the estimators*

| Estimator | $\hat{\sigma}_{11}$ | $\hat{\sigma}_{12}$ | $\hat{\sigma}_{13}$ | $\hat{\sigma}_{22}$ | $\hat{\sigma}_{23}$ | $\hat{\sigma}_{33}$ | $\hat{\sigma}^2$ |
|---|---|---|---|---|---|---|---|
| MSE | 0.0055 | 0.0013 | 0.0013 | 0.0052 | 0.0015 | 0.0054 | 0.0029 |

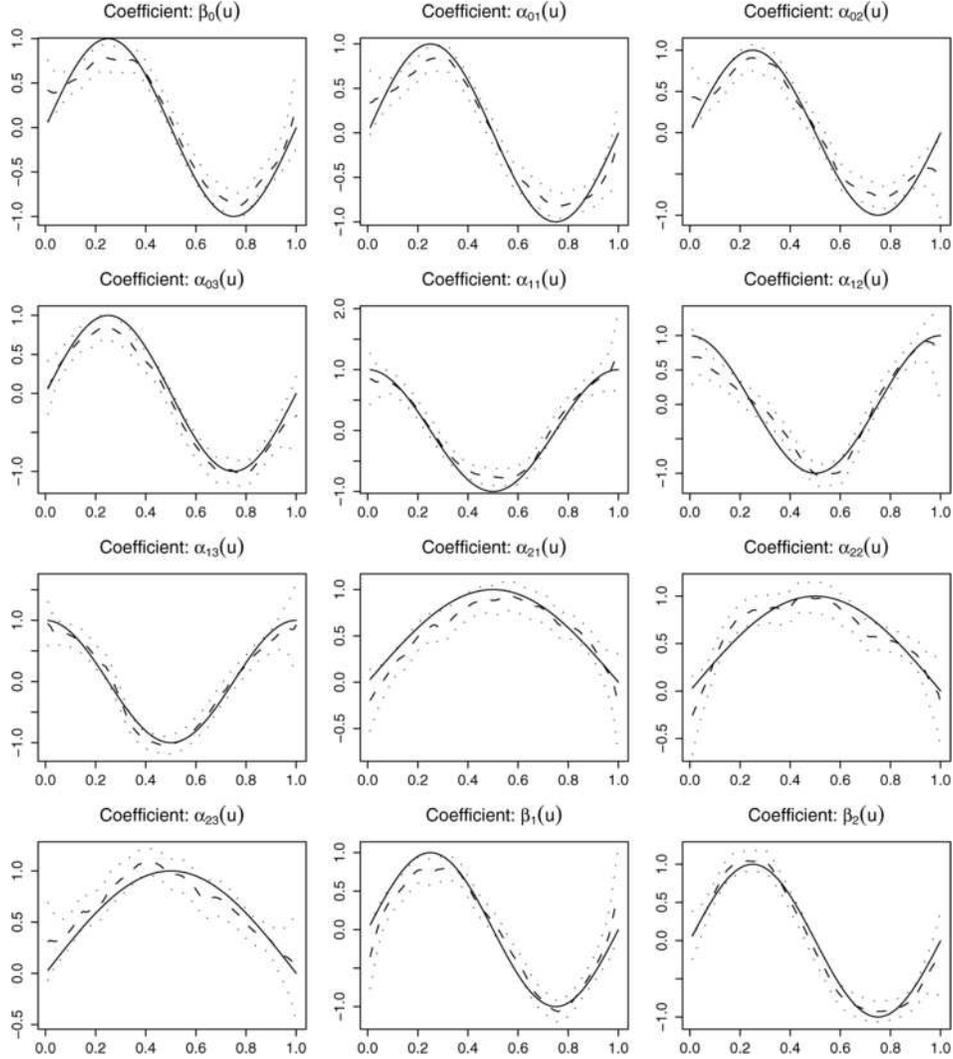

Fig. 1. *The solid lines are the true curves, the dashed lines are the estimators, and the dotted lines are* 95% *confidence bands.*

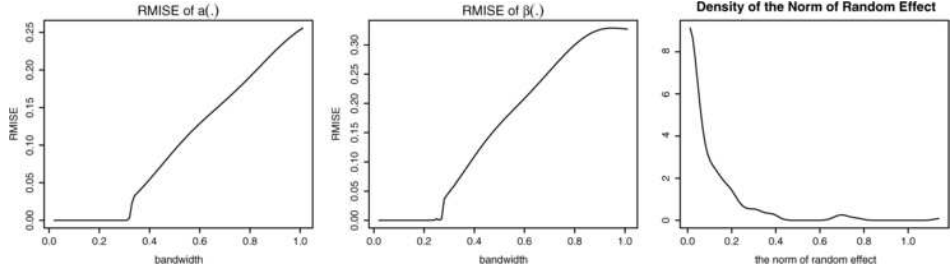

FIG. 2. *The left panel is the plot of the RMISEs of $\mathbf{a}(\cdot)$ against bandwidths, the middle panel is the plot of the RMISEs of $\boldsymbol{\beta}(\cdot)$, and the right panel is the density function of the estimated random effects in real data analysis.*

**6. Real data analysis.** The data we study come from the Bangladesh Demographic and Health Survey (BDHS) of 1996–1997 (Mitra et al. [22]), which is a cross-sectional, nationally representative survey of ever-married women aged between 10 and 49. The analysis is based on a sample of 8189 women nested within 296 primary sampling units, or clusters, with sample sizes ranging from 16 to 58. We allow the hierarchical structure of the data by fitting a two-level model, with women at level 1 nested within clusters at level 2. A further hierarchical level is the administrative division. Bangladesh is divided into six administrative divisions: Barisal, Chittagong, Dhaka, Kulna, Rajshahi and Sylhet. Effects at this level are represented in the model by fixed effects, since there are only six divisions.

The dependent variable $y_{ij}$ is the duration in months between the first birth and the second birth for the $j$th woman in the $i$th cluster. We consider several covariates that are commonly found to be associated with fertility in Bangladesh. The selected individual-level categorical covariates include the woman's level of education (none coded by 0 and primary or secondary plus coded by 1) denoted by $x_{ij1}$, religion (Hindu coded by 1 and Muslim or other coded by 0) denoted by $x_{ij2}$, first child (boy coded by 1) denoted by $x_{ij3}$. $X_{ij} = (x_{ij1}, x_{ij2}, x_{ij3})^T$. Another individual-level covariate is the year of marriage ($U_{ij}$). We also consider two cluster-level variables, administrative division and type of region of residence (rural coded by 1 and urban coded by 0). We take urban as the reference, and the differences between urban and rural clusters are modeled by dummy variables $z_{i1}$. We take Barisal as the reference, and the differences among the six administrative divisions are modeled by a set of dummy variables $z_{il}$, $l = 2, \ldots, 6$.

The proposed model (1.4) is used to fit the data, and the proposed estimation procedure is employed to estimate $\boldsymbol{\alpha}_j(\cdot) = (\alpha_{j1}(\cdot), \alpha_{j2}(\cdot), \alpha_{j3}(\cdot))^T$, $j = 0, \ldots, 6$ and $\boldsymbol{\beta}(\cdot) = (\beta_1(\cdot), \ldots, \beta_6(\cdot))^T$. The kernel involved in the estimation is taken to be Epanechnikov kernel, and the bandwidth is chosen to be 35% of the range of $U_{ij}$.

TABLE 3
*The P-values of the coefficients being constant*

| Coefficient | $\beta_0(\cdot)$ | $\alpha_{01}(\cdot)$ | $\alpha_{02}(\cdot)$ | $\alpha_{03}(\cdot)$ | $\alpha_{11}(\cdot)$ | $\alpha_{12}(\cdot)$ |
|---|---|---|---|---|---|---|
| P-value | 0.102 | 0.005 | 0.008 | 0.248 | 0.094 | 0.217 |
| Coefficient | $\alpha_{13}(\cdot)$ | $\alpha_{21}(\cdot)$ | $\alpha_{22}(\cdot)$ | $\alpha_{23}(\cdot)$ | $\alpha_{31}(\cdot)$ | $\alpha_{32}(\cdot)$ |
| P-value | 0.082 | 0.074 | 0.213 | 0.060 | 0.618 | 0.038 |
| Coefficient | $\alpha_{33}(\cdot)$ | $\alpha_{41}(\cdot)$ | $\alpha_{42}(\cdot)$ | $\alpha_{43}(\cdot)$ | $\alpha_{51}(\cdot)$ | $\alpha_{52}(\cdot)$ |
| P-value | 0.283 | 0.235 | 0.031 | 0.135 | 0.263 | 0.020 |
| Coefficient | $\alpha_{53}(\cdot)$ | $\alpha_{61}(\cdot)$ | $\alpha_{62}(\cdot)$ | $\alpha_{63}(\cdot)$ | $\beta_1(\cdot)$ | $\beta_2(\cdot)$ |
| P-value | 0.052 | 0.052 | 0.148 | 0.269 | 0.372 | 0.174 |
| Coefficient | $\beta_3(\cdot)$ | $\beta_4(\cdot)$ | $\beta_5(\cdot)$ | $\beta_6(\cdot)$ | | |
| P-value | 0.069 | 0.035 | 0.023 | 0.073 | | |

First, we examine how strong the random effect on each cluster is. To this end, for each cluster, we estimate the random effect of this cluster. The density function of the norm of the random effects of all clusters is presented in Figure 2. It is easy to see that the random effects are close to zero. This indicates that the within-cluster correlation has been mainly accounted by the deterministic cluster effect in (1.4).

As we mentioned before, if a coefficient is treated as a function when it is constant, we would pay a price on the variances of the resulting estimators. Thus, for each coefficient in the model, the proposed hypothesis test is employed to test whether it is constant or not. The P-value for each coefficient is depicted in Table 3. Table 3 shows that $\alpha_{01}(\cdot)$, $\alpha_{02}(\cdot)$, $\alpha_{32}(\cdot)$, $\alpha_{42}(\cdot)$, $\alpha_{52}(\cdot)$, $\beta_4(\cdot)$ and $\beta_5(\cdot)$ are nonconstant, and the others are constant. From now on, we use $\alpha_{ij}$ to denote constant function $\alpha_{ij}(\cdot)$ and $\beta_i$ for constant function $\beta_i(\cdot)$.

The former indicates the presentness of their nonlinear interactions with the year of marriage, while the latter shows no such interactions to be present.

The proposed estimation procedure is used to estimate the constant and functional coefficients. The estimated constant coefficients and their standard errors computed by the leave-one-cluster-out Jackknife are presented in Table 4, and the estimated functional coefficients together with their 95% confidence bands are presented in Figure 3.

It is visible, from Table 4, that some constant coefficients are not significantly apart from zero. This has some practical indications. For example, $\alpha_{41}$ (its estimated value is $-0.014$ with standard error 0.016) not significantly apart from zero indicates that the impact of eduction in the division of Kulna is not significantly different from that in Barisal.

As explained in the section of introduction, the proposed modeling and estimation methods mainly serve for the inference for a particular cluster.

We are now going to use a cluster in the rural area in Chittagong to illustrate how the proposed method works.

Based on the model (1.4) and the estimators of the unknown coefficients involved, we have the impact of education on the second birth interval in a rural area in Chittagong

$$\hat{a}_1(U) = \hat{\alpha}_{01}(U) + \hat{\alpha}_{11} + \hat{\alpha}_{21}.$$

Similarly, we can get the impact of Hindu

$$\hat{a}_2(U) = \hat{\alpha}_{02}(U) + \hat{\alpha}_{12} + \hat{\alpha}_{22}$$

and the impact of first child in this area

$$\hat{a}_3(U) = \hat{\alpha}_{03} + \hat{\alpha}_{13} + \hat{\alpha}_{23}.$$

The functional coefficients $\hat{a}_1(\cdot)$ and $\hat{a}_2(\cdot)$ are presented in Figure 3. The second one in the bottom panel in Figure 3 is $\hat{a}_1(\cdot)$, which shows that, in rural area of Chittagong, even the educated women still have a shorter second birth interval than the uneducated women in urban area in Barisal before 1970. This indicates that administrative division and the type of region of residence play a very important role in fertility behavior in Bangladesh before 1970. It is also noticeable that the second birth intervals of the educated women are getting longer.

From the third one in the bottom panel in Figure 3, which is $\hat{a}_2(\cdot)$, we can see that, even in rural areas in Chittagong, Hindus still have longer second birth intervals than Muslims even in the urban area in Barisal after 1970. This suggests that religion plays an important role in fertility behavior in Bangladesh.

The impact $\hat{a}_3(\cdot)$ of the first child being a boy in Rural in Chittagong does not vary with time. It is 0.128, which suggests that, even in the rural area of Chittagong, the women with a first child being a boy still have longer

TABLE 4
*Estimated constant coefficients*

| Coefficient | $\beta_0$ | $\alpha_{03}$ | $\alpha_{11}$ | $\alpha_{12}$ | $\alpha_{13}$ | $\alpha_{21}$ | $\alpha_{22}$ |
|---|---|---|---|---|---|---|---|
| Estimate | 3.530 | 0.016 | $-0.037$ | 0.105 | 0.011 | 0.036 | 0.134 |
| SE | 0.004 | 0.005 | 0.010 | 0.004 | 0.003 | 0.013 | 0.006 |
| Coefficient | $\alpha_{23}$ | $\alpha_{31}$ | $\alpha_{33}$ | $\alpha_{41}$ | $\alpha_{43}$ | $\alpha_{51}$ | $\alpha_{53}$ |
| Estimate | 0.101 | $-0.060$ | 0.042 | $-0.014$ | 0.111 | $-0.091$ | 0.104 |
| SE | 0.005 | 0.006 | 0.005 | 0.016 | 0.007 | 0.009 | 0.005 |
| Coefficient | $\alpha_{61}$ | $\alpha_{62}$ | $\alpha_{63}$ | $\beta_1$ | $\beta_2$ | $\beta_3$ | $\beta_6$ |
| Estimate | $-0.048$ | 0.171 | $-0.005$ | $-0.043$ | $-0.103$ | $-0.022$ | $-0.093$ |
| SE | 4.149 | 0.010 | 0.007 | 0.002 | 0.004 | 0.004 | 0.006 |

Note: SE stands for standard error of the estimator.

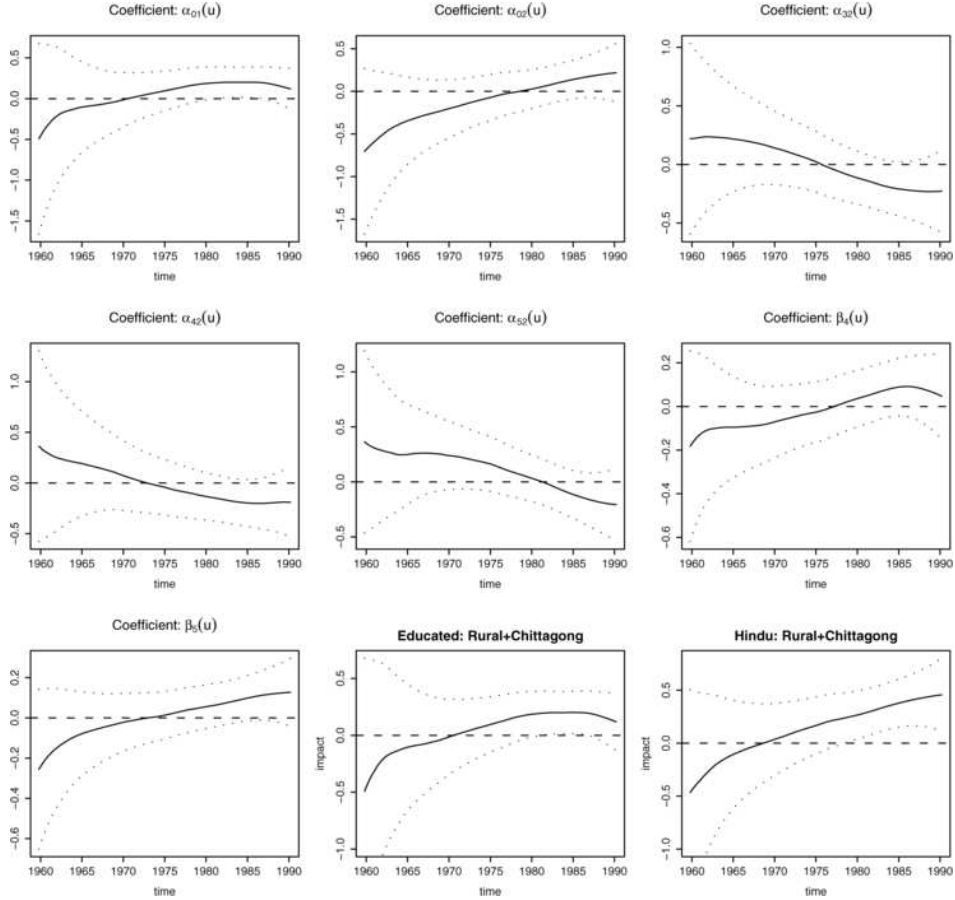

Fig. 3. *The solid lines are the estimated functional coefficients, and dotted lines are the 95% confidence bands.*

second birth intervals than the women with first child being a girl or dead even in the urban area of Barisal. A possible interpretation is that, like many developing countries, the Bangladesh culture always favors a boy.

## APPENDIX

In this section, we will prove the asymptotic distributions of the proposed estimators. For easy description, we write

$$\boldsymbol{\varepsilon} = (\boldsymbol{\varepsilon}_1^T, \ldots, \boldsymbol{\varepsilon}_m^T)^T, \qquad \mathbf{x} = \mathrm{diag}(\mathbf{x}_1, \ldots, \mathbf{x}_m), \qquad \mathbf{e} = (\mathbf{e}_1^T, \ldots, \mathbf{e}_m^T)^T,$$

$$\hat{\boldsymbol{\theta}}(u) = (\hat{\boldsymbol{\alpha}}_0(u)^T, \hat{\boldsymbol{\alpha}}_1(u)^T, \ldots, \hat{\boldsymbol{\alpha}}_q(u)^T, \hat{\boldsymbol{\beta}}(u)^T)^T,$$

$$\tau(u) = \{e_{(k+1),(q+1)}^T \otimes I_p\}[\sigma^2 \Lambda_1(u)^{-1} + \Upsilon_1(u) \Xi_1(u) \Upsilon_1(u)^T]\{e_{(k+1),(q+1)} \otimes I_p\}$$

and

$$H = \begin{pmatrix} 1 & 0 \\ 0 & h \end{pmatrix} \otimes I_s,$$

where $s = (q+1)p + q$ is defined in (2.2).

Moreover, for any function $g(u)$ on the interval $[a, b]$, define $\|g\|_\infty = \sup_{u \in [a,b]} |g(u)|$, and, for any matrix $A(u) = (a_{ij}(u))_{p \times p}$, set

$$\|A\|_\infty = \left( \sum_{i=1}^p \sum_{j=1}^p \|a_{ij}\|_\infty^2 \right)^{1/2}.$$

The following technical conditions are imposed to establish the asymptotic results:

(1) $E\varepsilon_{11}^4 < \infty, E\|\mathbf{e}_1\|^4 < \infty, Ez_j^4 < \infty, Ex_i^{2d} < \infty$ for some $d > 2$, where $\|\mathbf{e}_1\|^2 = \mathbf{e}_1^T \mathbf{e}_1$, $z_j$ denotes the $j$th element of $Z$ and $x_i$ denotes the $i$th element of $X$ for $j = 1, \ldots, q, i = 1, \ldots, p$;
(2) $\ddot{\alpha}_{kj}(\cdot)$ is continuous in a neighborhood of $u$, for $k = 0, \ldots, q, j = 1, \ldots, p$, where $\ddot{\alpha}_{kj}(\cdot)$ is the $j$th element of $\ddot{\boldsymbol{\alpha}}_k(\cdot)$, and assume $\ddot{\alpha}_{kj}(u) \neq 0$; similarly, $\ddot{\beta}_l(\cdot)$ is continuous in a neighborhood of $u$, for $l = 1, \ldots, q$, where $\ddot{\beta}_l(\cdot)$ is the $l$th element of $\boldsymbol{\beta}(\cdot)$, and $\ddot{\beta}_l(u) \neq 0$;
(3) The marginal density $f(\cdot)$ of $U$ has a continuous derivative in some neighborhood of $u$, and $f(u) \neq 0$;
(4) Each element of $\Omega(u)$ and $\Xi_1(u)$ are continuous in the neighborhood of $u$, and $\Omega(u)$ is positive definite at the point $u$;
(5) The function $K(t)$ is a symmetric density function with a compact support;
(6) $n_i, i = 1, \ldots, m$, are bounded. $n \to \infty$, $h \to 0$ and $nh^2 \to \infty$;
(7) $E\|(\mathbf{x}_i^T \mathbf{x}_i)^{-1}\|^{4+\delta} < \infty$ for some $\delta > 0$ and $i = 1, \ldots, m$, where $\|\cdot\|$ denotes Frobenium norm of a matrix;
(8) $E[X_{il}^T \Sigma X_{ir} \Gamma_{il} \Gamma_{ir}^T | U_{il} = u, U_{ir} = v]$ is continuous in the neighborhood of $u$ and $v$ respectively for $r \neq l, r, l = 1, \ldots, n_i, i = 1, \ldots, m$;
(9) $\Omega(u)$ is bounded and has a continuous derivative on $[a, b]$, $\tau(u)$ has a bounded derivative on $[a, b]$, $0 < \|\tau(u)\|_\infty < \infty$ and the first derivative of $K(t)$ has a finite number of sign changes over its support.

Note that $\Omega(u)$ is automatically positive semidefinite at the point $u$, so the second part of condition (4) is easily satisfied.

To obtain the proof of the theorems, the following lemmas are required.

LEMMA 1. *Let $\{U_{ij}\}$ be i.i.d. random variables, $\{\xi_{ij}\}$ be identically distributed and $\xi_{ij}$ be independent of $\xi_{lk}$ for $i \neq l$. Assume that the marginal density $f(\cdot)$ of $U$ has a continuous derivative in some neighborhood of $u$, $E(\xi_{11}|U_{11} = u)$ is continuous in the neighborhood of $u$ and $E(\xi_{11}^2) < \infty$. Let* $K(\cdot)$ be a bounded positive function with a bounded support. Then, when $nh^2 \to \infty$, for $\lambda = 0, 1, 2, \ldots$,

$$n^{-1} \sum_{i=1}^{m} \sum_{j=1}^{n_i} \xi_{ij} (h^{-1}(U_{ij} - u))^\lambda K_h(U_{ij} - u)$$
$$= \mu_\lambda f(u) E(\xi_{11} | U_{11} = u)(1 + o_P(1)).$$

PROOF. Let $S_{n,\lambda} = n^{-1} \sum_{i=1}^{m} \sum_{j=1}^{n_i} \xi_{ij}(h^{-1}(U_{ij} - u))^\lambda K_h(U_{ij} - u)$ and $g(v) = E(\xi_{11} | U_{11} = v)$. Then,

$$ES_{n,\lambda} = n^{-1} \sum_{i=1}^{m} \sum_{j=1}^{n_i} E\{E[\xi_{ij}(h^{-1}(U_{ij} - u))^\lambda K_h(U_{ij} - u) | U_{ij}]\}$$
$$= \int g(v)(h^{-1}(v - u))^\lambda K_h(v - u) f(v) \, dv = \mu_\lambda f(u) g(u)(1 + o(1))$$

by continuity of both the density function $f(\cdot)$ and the conditional expectation function $g(\cdot)$ in the neighborhood of $u$. Moreover, as $K(\cdot)$ is a bounded function with a bounded support, $|u^\lambda K(u)|$ is bounded. Then, by the Jensen inequality, it follows that

$$\text{var}(S_{n,\lambda}) \leq ES_{n,\lambda}^2 \leq n^{-2} \sum_{i=1}^{m} n_i \sum_{j=1}^{n_i} E[\xi_{ij}^2(h^{-1}(U_{ij} - u))^{2\lambda} K_h^2(U_{ij} - u)]$$
$$= O((nh^2)^{-1}).$$

Therefore,

$$S_{n,\lambda} = ES_{n,\lambda} + O_P(\sqrt{\text{var}(S_{n,\lambda})}) = \mu_\lambda f(u) E(\xi_{11} | U_{11} = u)(1 + o_P(1)). \quad \Box$$

LEMMA 2. *Let $\{\xi_n, n \geq 1\}$ be an independent sequence, with $E\xi_n = 0$ and*

$$\lambda E|\xi_j|^3 \exp(\lambda|\xi_j|) \leq E\xi_j^2$$

*for any $j \geq 1$ and some $\lambda > 0$. If $\sum_{i=1}^{n} E\xi_i^2 \to \infty$, then, on a possibly enlarged probability space, there exists a sequence of independent random variables $\{\eta_n, n \geq 1\}, \eta_n \sim N(0, \text{var}(\xi_n))$ such that*

$$\left| \sum_{i=1}^{n} \xi_i - \sum_{i=1}^{n} \eta_i \right| \leq \frac{1}{\lambda C} \log\left( \sum_{i=1}^{n} \text{var}(\xi_i) \right),$$

*where $C$ is a positive constant.*

See Lin and Lu [20], page 129, Theorem 2.6.3.

PROOF OF THEOREM 1. It can be shown that

$$\sqrt{nhf(u)}(\hat{\boldsymbol{\alpha}}_k(u) - \boldsymbol{\alpha}_k(u))$$
$$= \sqrt{nhf(u)}(e_{(k+1),(q+1)}^T \otimes I_p, \mathbf{0}_{p \times (q+s)})(\mathbf{X}^T W \mathbf{X})^{-1} \mathbf{X}^T W \boldsymbol{\varepsilon}$$
$$+ \sqrt{nhf(u)}(e_{(k+1),(q+1)}^T \otimes I_p, \mathbf{0}_{p \times (q+s)})(\mathbf{X}^T W \mathbf{X})^{-1} \mathbf{X}^T W \mathbf{xe}$$
$$+ \sqrt{nhf(u)}\{(e_{(k+1),(q+1)}^T \otimes I_p, \mathbf{0}_{p \times (q+s)})(\mathbf{X}^T W \mathbf{X})^{-1} \mathbf{X}^T W E(Y|\mathcal{D})$$
$$- \boldsymbol{\alpha}_k(u)\}$$
$$\equiv L_{n1} + L_{n2} + L_{n3},$$

as

$$L_{n1} = n(e_{(k+1),(q+1)}^T \otimes I_p, \mathbf{0}_{p \times (q+s)})(\mathbf{X}^T W \mathbf{X})^{-1} H \sqrt{n^{-1} h f(u)}\{H^{-1} \mathbf{X}^T W \boldsymbol{\varepsilon}\}$$

with

$$\sqrt{n^{-1} h f(u)} E\{H^{-1} \mathbf{X}^T W \boldsymbol{\varepsilon}\} = \mathbf{0}_{2s \times 1}$$

and

$$n^{-1} h f(u) \operatorname{cov}\{H^{-1} \mathbf{X}^T W \boldsymbol{\varepsilon}\} = \sigma^2 f(u) n^{-1} h E\{H^{-1} X^T W^2 X H^{-1}\}$$
$$= \sigma^2 f^2(u) \begin{pmatrix} \nu_0 & 0 \\ 0 & \nu_2 \end{pmatrix} \otimes \Omega(u)(1 + o(1)).$$

Moreover, it follows from Lemma 1 that

$$(A.1) \quad \frac{1}{n}(\mathbf{X}^T W \mathbf{X}) = f(u) H \left\{ \begin{pmatrix} 1 & 0 \\ 0 & \mu_2 \end{pmatrix} \otimes \Omega(u) \right\} H(1 + o_P(1)).$$

Hence,

$$n(e_{(k+1),(q+1)}^T \otimes I_p, \mathbf{0}_{p \times (q+s)})(\mathbf{X}^T W \mathbf{X})^{-1} H$$
$$= \frac{1}{f(u)}((e_{(k+1),(q+1)}^T \otimes I_p, \mathbf{0}_{p \times q})\Omega(u)^{-1}, \mathbf{0}_{p \times s})(1 + o_P(1)).$$

By conditions (1), (5) and (6), Lindeberg–Feller Theorem, Slutsky's theorem and inverse of block matrix, it follows that

$$L_{n1} \xrightarrow{D} N_p(\mathbf{0}_{p \times 1}, \nu_0 \sigma^2 \{e_{(k+1),(q+1)}^T \otimes I_p\} \Lambda_1(u)^{-1} \{e_{(k+1),(q+1)} \otimes I_p\}).$$

Similarly, it can be shown that

$$L_{n2} \xrightarrow{D} N_p(\mathbf{0}_{p \times 1}, \nu_0 \{e_{(k+1),(q+1)}^T \otimes I_p\} \Theta_1(u) \{e_{(k+1),(q+1)} \otimes I_p\}).$$

Since $L_{n1}$ and $L_{n2}$ are independent, $(L_{n1} + L_{n2})$ has the asymptotic distribution

$$N_p(\mathbf{0}_{p \times 1}, \nu_0 \{e_{(k+1),(q+1)}^T \otimes I_p\}[\sigma^2 \Lambda_1(u)^{-1} + \Theta_1(u)]\{e_{(k+1),(q+1)} \otimes I_p\}).$$

By similar arguments as establishment of Lemma 1 and condition (2), we find that

$$\mathbf{X}^T W E(Y|\mathcal{D}) - \mathbf{X}^T W \mathbf{X} \begin{pmatrix} \boldsymbol{\theta}(u) \\ \dot{\boldsymbol{\theta}}(u) \end{pmatrix}$$

(A.2)
$$= \begin{pmatrix} \frac{1}{2} n h^2 \mu_2 f(u) \Omega(u) \\ \mathbf{0}_{s \times s} \end{pmatrix} \ddot{\boldsymbol{\theta}}(u)(1 + o_P(1)).$$

This, together with (A.1), we get that

$$L_{n3} = \frac{\mu_2 \sqrt{nh^5 f(u)}}{2} \ddot{\boldsymbol{\alpha}}_k(u)(1 + o_P(1)).$$

Therefore, when $nh^5$ is bounded, we obtain that

$$\sqrt{nhf(u)} \left\{ \hat{\boldsymbol{\alpha}}_k(u) - \boldsymbol{\alpha}_k(u) - h^2 \frac{\mu_2}{2} \ddot{\boldsymbol{\alpha}}_k(u) \right\}$$

$$\xrightarrow{D} N_p(\mathbf{0}_{p \times 1},$$
$$\nu_0 \{e_{(k+1),(q+1)}^T \otimes I_p\} [\sigma^2 \Lambda_1(u)^{-1} + \Theta_1(u)] \{e_{(k+1),(q+1)} \otimes I_p\}).$$

□

PROOF OF THEOREM 2. It follows immediately from the proof of Theorem 1. □

PROOF OF THEOREM 3. Let $\Delta \mathbf{r}_i = \hat{\mathbf{r}}_i - \mathbf{r}_i$, $Q_i = I_{n_i} - P_i$ and $\tilde{Q}_i$ be a diagonal matrix generated from the diagonal elements of $Q_i$. Now,

$$\hat{\sigma}^2 = \frac{1}{n - mp} \sum_{i=1}^m \boldsymbol{\varepsilon}_i^T (Q_i - \tilde{Q}_i) \boldsymbol{\varepsilon}_i + \frac{1}{n - mp} \sum_{i=1}^m \boldsymbol{\varepsilon}_i^T \tilde{Q}_i \boldsymbol{\varepsilon}_i$$

$$+ \frac{1}{n - mp} \sum_{i=1}^m E[\Delta \mathbf{r}_i^T | \mathcal{D}] Q_i E[\Delta \mathbf{r}_i | \mathcal{D}]$$

$$+ \frac{1}{n - mp} \sum_{i=1}^m \{\Delta \mathbf{r}_i - E[\Delta \mathbf{r}_i | \mathcal{D}]\}^T Q_i \{\Delta \mathbf{r}_i - E[\Delta \mathbf{r}_i | \mathcal{D}]\}$$

(A.3)
$$+ \frac{2}{n - mp} \sum_{i=1}^m \{\Delta \mathbf{r}_i - E[\Delta \mathbf{r}_i | \mathcal{D}]\}^T Q_i E[\Delta \mathbf{r}_i | \mathcal{D}]$$

$$+ \frac{2}{n - mp} \sum_{i=1}^m \boldsymbol{\varepsilon}_i^T Q_i E[\Delta \mathbf{r}_i | \mathcal{D}]$$

$$+ \frac{2}{n - mp} \sum_{i=1}^m \boldsymbol{\varepsilon}_i^T Q_i \{\Delta \mathbf{r}_i - E[\Delta \mathbf{r}_i | \mathcal{D}]\}$$

$$\equiv J_{n1} + J_{n2} + J_{n3} + J_{n4} + J_{n5} + J_{n6} + J_{n7}.$$

As $Q_i$ is an idempotent matrix and all the diagonal components of $Q_i - \tilde{Q}_i$ are equal to zero, by straightforward calculation, it follows that

$$E(J_{n1}|\mathcal{D}) = \frac{\sigma^2}{n-mp} \sum_{i=1}^{m} \text{tr}(Q_i - \tilde{Q}_i) = 0,$$

$$E(J_{n2}|\mathcal{D}) = \frac{\sigma^2}{n-mp} \sum_{i=1}^{m} \text{tr}(\tilde{Q}_i) = \sigma^2,$$

$$\text{cov}(J_{n1}, J_{n2}|\mathcal{D}) = E(J_{n_1} J_{n_2}|\mathcal{D}) = \frac{2\sigma^4}{(n-mp)^2} \sum_{i=1}^{m} \text{tr}((Q_i - \tilde{Q}_i)\tilde{Q}_i) = 0$$

and

$$\text{var}(J_{n1}) = E\{E(J_{n1}^2|\mathcal{D})\} = \frac{2\sigma^4}{n-mp} E\left\{\frac{mp - \sum_{i=1}^{m}\sum_{j=1}^{n_i}[X_{ij}^T(\mathbf{x}_i^T\mathbf{x}_i)^{-1}X_{ij}]^2}{n-mp}\right\},$$

$$\text{var}(J_{n2}) = E\{E(J_{n2}^2|\mathcal{D})\} - \sigma^4 = \frac{\text{var}(\varepsilon_{11}^2)}{(n-mp)^2} E\sum_{i=1}^{m}\sum_{j=1}^{n_i}[1 - X_{ij}^T(\mathbf{x}_i^T\mathbf{x}_i)^{-1}X_{ij}]^2.$$

As $J_{n1} = (n-mp)^{-1}\sum_{i=1}^{m}\{\sum_{l=1}^{n_i}\sum_{r=1,r\neq l}^{n_i} X_{il}^T(\mathbf{x}_i^T\mathbf{x}_i)^{-1}X_{ir}\varepsilon_{il}\varepsilon_{ir}\}$ is a sum of independent variables, by Lindeberg–Feller Theorem, it follows that

$$n^{1/2} J_{n1} \xrightarrow{D} N(0, 2\sigma^4 c_1(c_1 - 1 - \gamma)).$$

Similarly,

$$n^{1/2}(J_{n2} - \sigma^2) \xrightarrow{D} N(0, \text{var}(\varepsilon_{11}^2)c_1(2 - c_1 + \gamma)).$$

Since the two terms are uncorrelated, we have that

(A.4)
$$n^{1/2}(J_{n1} + J_{n2} - \sigma^2)$$
$$\xrightarrow{D} N(0, 2\sigma^4 c_1(c_1 - 1 - \gamma) + \text{var}(\varepsilon_{11}^2)c_1(2 - c_1 + \gamma)).$$

In the following, we will show that the remaining parts from $J_{n3}$ to $J_{n7}$ in (A.3) satisfy $n^{1/2}J_{nl} = o_P(1), l = 3, \ldots, 7$.

By the conditional bias of $\hat{\boldsymbol{\theta}}(u)$ and law of large numbers, it follows from $nh^8 \to 0$ that

(A.5)
$$n^{1/2} J_{n3} = o_P(1).$$

Since $0 \leq \boldsymbol{\phi}^T Q_i \boldsymbol{\phi} \leq \boldsymbol{\phi}^T \boldsymbol{\phi}$ for any $i$ and $n_i$ dimensional vector $\boldsymbol{\phi}$, we have

$$E\{|J_{n4}||\mathcal{D}\} \leq (n-mp)^{-1} \sum_{i=1}^{m}\sum_{j=1}^{n_i} \Gamma_{ij}^T \text{cov}(\hat{\boldsymbol{\theta}}(U_{ij})|\mathcal{D})\Gamma_{ij} = O_P((nh)^{-1}).$$

By the conditional bias and covariance matrix of $\hat{\boldsymbol{\theta}}(u)$, it follows that

$$E\{|J_{n5}||\mathcal{D}\}$$
$$= \frac{h^2\mu_2(1+o_P(1))}{(n-mp)}$$
$$\times \left\{ \sum_{i=1}^{m}\sum_{k=1}^{n_i-p}\sum_{j=1}^{n_i}\sum_{r=1}^{n_i} |Q_{ikj}Q_{ikr}\Gamma_{ir}^T\ddot{\boldsymbol{\theta}}(U_{ir})| \right.$$
$$\left. \times E[|\Gamma_{ij}^T(\hat{\boldsymbol{\theta}}(U_{ij}) - E[\hat{\boldsymbol{\theta}}(U_{ij})|\mathcal{D}])||\mathcal{D}] \right\}$$
$$\leq \frac{h^2\mu_2(1+o_P(1))}{(n-mp)}$$
$$\times \sum_{i=1}^{m}(n_i-p)\left\{ \sum_{r=1}^{n_i}(\Gamma_{ir}^T\ddot{\boldsymbol{\theta}}(U_{ir}))^2 \sum_{j=1}^{n_i}\Gamma_{ij}^T\operatorname{cov}(\hat{\boldsymbol{\theta}}(U_{ij})|\mathcal{D})\Gamma_{ij} \right\}^{1/2}$$
$$= O_p((n^{-1}h^3)^{1/2}).$$

where

$$\sum_{l=1}^{n_i} Q_{irl}Q_{ivl} = \delta_{rv} = \begin{cases} 1, & r = v, \\ 0, & r \neq v, \end{cases}$$

and

$$E(J_{n6}^2|\mathcal{D}) \leq \frac{h^4\mu_2^2\sigma^2}{(n-mp)^2} \sum_{i=1}^{m}\sum_{j=1}^{n_i}[\Gamma_{ij}^T\ddot{\boldsymbol{\theta}}(U_{ij})]^2(1+o_P(1)).$$

By a similar expression as (A.1) and straightforward calculation, we have

$$\Gamma_{ij}^T[\hat{\boldsymbol{\theta}}(U_{ij}) - E(\hat{\boldsymbol{\theta}}(U_{ij})|\mathcal{D})]\varepsilon_{ir}$$
$$= \frac{1}{nf(U_{ij})}\Gamma_{ij}^T\Omega(U_{ij})^{-1}\left(\sum_{t=1}^{m}\sum_{l=1}^{n_t}\Gamma_{tl}\varepsilon_{tl}K_h(U_{tl}-U_{ij})\right)\varepsilon_{ir}(1+o_P(1)).$$

Moreover, by boundness of the kernel function and independence of the random errors, it can be shown that $E\{|J_{n7}||\mathcal{D}\} = O_P((nh)^{-1})$.

Therefore, using Markov inequality, when $h \to 0, nh^2 \to \infty$ we get

(A.6) $$n^{1/2}J_{nl} = o_P(1), \quad l = 4,\ldots,7.$$

Combing the results from (A.3) to (A.6), we have

$$n^{1/2}\{\hat{\sigma}^2 - \sigma^2\} \xrightarrow{D} N(0, 2\sigma^4 c_1(c_1-1-\gamma) + \operatorname{var}(\varepsilon_{11}^2)c_1(2-c_1+\gamma)). \quad \square$$

PROOF OF THEOREM 4. Using standard arguments as in the proof of Theorem 3 and the law of large numbers, when $nh^2 \to \infty$, the conditional bias of $\hat{\Sigma}$ is

$$\text{bias}\{\text{vec}(\hat{\Sigma})|\mathcal{D}\} = O_P(h^4) + o_P(n^{-1/2}),$$

and, by straightforward but tedious calculation and Lindeberg–Feller Theorem, when $nh^8 \to 0$, it follows that

$$n^{1/2}\,\text{vec}(\hat{\Sigma} - \Sigma) \xrightarrow{D} N_{p^2}(\mathbf{0}_{p^2 \times 1}, (1/c_2 + p)\Delta)$$

where

$$\Delta = E\{(\mathbf{e}_1\mathbf{e}_1^T) \otimes (\mathbf{e}_1\mathbf{e}_1^T)\} - \text{vec}(\Sigma)\text{vec}(\Sigma)^T + 2\sigma^4 \Delta_2$$
$$+ \sigma^2\{\Sigma \otimes \Delta_1 + \Delta_1 \otimes \Sigma + \Delta_3\}$$
$$+ [\text{var}(\varepsilon_{11}^2) - 2\sigma^4]\{\Delta_4 - c_2[2\,\text{vec}(\Delta_1)\,\text{vec}(\Delta_1)^T$$
$$- \text{vec}(\Delta_1)\text{vec}(\Delta_2)^T - \text{vec}(\Delta_2)\text{vec}(\Delta_1)^T]\}$$
$$+ \{2\sigma^4 c_1(c_1 - 1 - \gamma) + \text{var}(\varepsilon_{11}^2)c_1(2 - c_1 + \gamma)\}\text{vec}(\Delta_1)\text{vec}(\Delta_1)^T.$$

Therefore, we have

$$n^{1/2}\,\text{vech}(\hat{\Sigma} - \Sigma)$$
$$\xrightarrow{D} N_{p(p+1)/2}(\mathbf{0}_{(p(p+1)/2) \times 1}, (1/c_2 + p)(R_p^T R_p)^{-1} R_p^T \Delta R_p (R_p^T R_p)^{-1}).$$

□

PROOF OF THEOREM 5. By simple calculation, it can be seen that

$$\hat{C}_{kj} - C_{kj}$$
$$= \frac{1}{n}\sum_{i=1}^m \sum_{l=1}^{n_i} e_{kp+j,2s}^T (\mathbf{X}_{(il)}^T W_{(il)} \mathbf{X}_{(il)})^{-1} \mathbf{X}_{(il)}^T W_{(il)} \boldsymbol{\varepsilon}$$
$$+ \frac{1}{n}\sum_{i=1}^m \sum_{l=1}^{n_i} e_{kp+j,2s}^T (\mathbf{X}_{(il)}^T W_{(il)} \mathbf{X}_{(il)})^{-1} \mathbf{X}_{(il)}^T W_{(il)} \mathbf{xe}$$
$$+ \frac{1}{n}\sum_{i=1}^m \sum_{l=1}^{n_i} e_{kp+j,2s}^T (\mathbf{X}_{(il)}^T W_{(il)} \mathbf{X}_{(il)})^{-1} \mathbf{X}_{(il)}^T W_{(il)}$$
$$\times \left\{ E(Y|\mathcal{D}) - \mathbf{X}_{(il)} \begin{pmatrix} \boldsymbol{\theta}(U_{il}) \\ \dot{\boldsymbol{\theta}}(U_{il}) \end{pmatrix} \right\}$$
$$\equiv T_{n1} + T_{n2} + T_{n3}.$$

First, we obtain that $E(T_{n1}) = 0$, and, for any $j$, let

$$F_{il}^T = e_{kp+j,2s}^T (\mathbf{X}_{(il)}^T W_{(il)} \mathbf{X}_{(il)})^{-1} \mathbf{X}_{(il)}^T W_{(il)}, \qquad l = 1, \ldots, n_i, i = 1, \ldots, m,$$

and $\mathbf{F} = (F_{11}, \ldots, F_{1n_1}, \ldots, F_{m1}, \ldots, F_{mn_m})^T$. Then, we have

$$\begin{aligned}
\operatorname{var}\{T_{n1}|\mathcal{D}\} &= \frac{\sigma^2}{n^2}\mathbf{1}_n^T\mathbf{F}\mathbf{F}^T\mathbf{1}_n \\
&= \frac{\sigma^2}{n^2}\sum_{i=1}^{m}\sum_{l=1}^{n_i}\sum_{r=1}^{m}\sum_{v=1}^{n_r} e_{kp+j,2s}^T(\mathbf{X}_{(il)}^T W_{(il)} \mathbf{X}_{(il)})^{-1} \\
&\qquad \times \mathbf{X}_{(il)}^T W_{(il)} W_{(rv)} \mathbf{X}_{(rv)} (\mathbf{X}_{(rv)}^T W_{(rv)} \mathbf{X}_{(rv)})^{-1} \\
&\qquad \times e_{kp+j,2s} \\
&= \frac{\sigma^2}{n} e_{kp+j,((q+1)p)}^T E\{\Lambda_1(U)^{-1}\} e_{kp+j,((q+1)p)} (1+o_P(1)),
\end{aligned}$$

by Lemma 1, straightforward but tedious calculation and the law of large numbers. Therefore, using conditions (1), (5) and (6) and the Lindeberg–Feller Theorem, it follows that

$$\sqrt{n}T_{n1} \xrightarrow{D} N(0, \sigma^2 e_{kp+j,((q+1)p)}^T E\{\Lambda_1(U)^{-1}\} e_{kp+j,((q+1)p)}).$$

By the same way, it can be shown that

$$\sqrt{n}T_{n2} \xrightarrow{D} N(0, e_{kp+j,((q+1)p)}^T \{E[\Theta_1(U)] + \Xi_2\} e_{kp+j,((q+1)p)}).$$

Moreover, combining the results similar to (A.1) and (A.2), we get

$$T_{n3} = h^2 \frac{\mu_2}{2}\left[\frac{1}{n}\sum_{i=1}^{m}\sum_{l=1}^{n_i} \ddot{\alpha}_{kj}(U_{il})\right](1+o_P(1)).$$

Therefore, by independence of $T_{n1}$ and $T_{n2}$, when $nh^4 \to 0$, we have

$$\sqrt{n}\{\hat{C}_{kj} - C_{kj}\}$$
$$\xrightarrow{D} N(0, e_{kp+j,((q+1)p)}^T \{\sigma^2 E[\Lambda_1(U)^{-1}] + E[\Theta_1(U)] + \Xi_2\} e_{kp+j,((q+1)p)}).$$

$\square$

PROOF OF THEOREM 6. Obviously,

$$\hat{\alpha}_{kj}(u) - \alpha_{kj}(u) - \operatorname{bias}\{\hat{\alpha}_{kj}(u)|\mathcal{D}\} = e_{kp+j,2s}^T(\mathbf{X}^T W \mathbf{X})^{-1}\mathbf{X}^T W(\boldsymbol{\varepsilon}+\mathbf{xe}) \equiv I_1(u).$$

First of all, we approximate the random matrix $I_1(u)$. As $f(\cdot)$ and $\Omega(\cdot)$ have continuous derivatives on $[a,b]$, by similar arguments as Lemma 1 and Neumann series expansion, it follows that

(A.7) $\quad nH(\mathbf{X}^T W \mathbf{X})^{-1} H = f(u)^{-1} S(u)^{-1} + O_P((nh^2)^{-1/2} + h)$

uniformly for $u \in [a,b]$ where $S(u) = \begin{pmatrix} \mu_0 & 0 \\ 0 & \mu_2 \end{pmatrix} \otimes \Omega(u)$.

By the asymptotic normality of $\sqrt{n^{-1}hf(u)}H^{-1}\mathbf{X}^T W(\boldsymbol{\varepsilon}+\mathbf{xe})$ in the proof of Theorem 1, we get that

$$\text{(A.8)} \qquad \left\|\frac{1}{n}H^{-1}\mathbf{X}^T W(\boldsymbol{\varepsilon}+\mathbf{xe})\right\|_\infty = O_P\left(\frac{1}{\sqrt{nh}}\right).$$

Therefore, using (A.7) and (A.8), it can be seen that

$$\text{(A.9)} \qquad \left\|I_1(u) - \frac{1}{n}e_{kp+j,2s}^T f(u)^{-1} S(u)^{-1} H^{-1}\mathbf{X}^T W(\boldsymbol{\varepsilon}+\mathbf{xe})\right\|_\infty$$
$$= O_P((nh^{3/2})^{-1} + (nh^{-1})^{-1/2}).$$

Next, we consider

$$I_2(u) \equiv \frac{1}{n}e_{kp+j,2s}^T f(u)^{-1} S(u)^{-1} H^{-1}\mathbf{X}^T W(\boldsymbol{\varepsilon}+\mathbf{xe})$$

$$= \sum_{i=1}^m \sum_{l=1}^{n_i} \frac{1}{nf(u)} e_{kp+j,s}^T \Omega(u)^{-1} \Gamma_{il} K_h(U_{il}-u)(\varepsilon_{il} + X_{il}^T \mathbf{e}_i).$$

Let $w_{il} = e_{kp+j,s}^T \Omega(u)^{-1} \Gamma_{il}$. Then,

$$\{nhf(u)\tau(u)^{-1}\nu_0^{-1}\}^{1/2} I_2(u)$$

$$\text{(A.10)} \qquad = \{nhf(u)\tau(u)\nu_0\}^{-1/2} \sum_{i=1}^m \sum_{l=1}^{n_i} K\left(\frac{U_{il}-u}{h}\right) w_{il}(\varepsilon_{il} + X_{il}^T \mathbf{e}_i)$$

$$\equiv \{nhf(u)\tau(u)\nu_0\}^{-1/2} I_3(u).$$

Divide interval $[a,b]$ into $n$ subintervals $J_r = [d_{r-1}, d_r), r = 1, 2, \ldots, n-1, J_n = [d_{n-1}, b]$ where $d_r = a + \frac{b-a}{n}r$. Define $\tilde{U}_{il} = d_r I(U_{il} \in J_r), r = 1, \ldots, n$, and it is obvious that $|U_{il} - \tilde{U}_{il}| = O(n^{-1})$. Then, by law of large numbers for random sequence $\{w_{il}(\varepsilon_{il} + X_{il}^T \mathbf{e}_i)\}$, it follows that

$$I_3(u) = \sum_{i=1}^m \sum_{l=1}^{n_i} \left[K\left(\frac{U_{il}-u}{h}\right) - K\left(\frac{\tilde{U}_{il}-u}{h}\right)\right] w_{il}(\varepsilon_{il} + X_{il}^T \mathbf{e}_i)$$

$$+ \sum_{i=1}^m \sum_{l=1}^{n_i} K\left(\frac{\tilde{U}_{il}-u}{h}\right) w_{il}(\varepsilon_{il} + X_{il}^T \mathbf{e}_i)$$

$$\text{(A.11)}$$

$$= O_P(h^{-1}) + \sum_{i=1}^m \sum_{l=1}^{n_i} K\left(\frac{\tilde{U}_{il}-u}{h}\right) w_{il}(\varepsilon_{il} + X_{il}^T \mathbf{e}_i)$$

$$\equiv O_P(h^{-1}) + I_4(u)$$

uniformly for $u \in [a,b]$. By the definition of $\tilde{U}_{il}$, we have that

$$I_4(u) = \sum_{r=1}^n K\left(\frac{d_r-u}{h}\right) \sum_{i=1}^m \sum_{l=1}^{n_i} I(U_{il} \in J_r) w_{il}(\varepsilon_{il} + X_{il}^T \mathbf{e}_i).$$

Let $\zeta_t = \sum_{r=1}^{t} \sum_{i=1}^{m} \sum_{l=1}^{n_i} I(U_{il} \in J_r) w_{il}(\varepsilon_{il} + X_{il}^T \mathbf{e}_i) = \sum_{i=1}^{m} \sum_{l=1}^{n_i} I(a \leq U_{il} < d_t) w_{il}(\varepsilon_{il} + X_{il}^T \mathbf{e}_i), \zeta_0 \equiv 0$. Then, by Lemma 2, for any $t = 1, \ldots, n$ and $u \in [a, b]$, we get

$$|\zeta_t - n^{1/2} W(G(d_t))| = O(n^{1/4} \log n) \quad \text{a.s.},$$

where $W(\cdot)$ is a Wiener process and

$$G(c) = \int_a^c \sigma^2 [E(w_{11}^2 | U_{11} = v) + E(X_{11}^T \Sigma X_{11} w_{11}^2 | U_{11} = v)] f(v) \, dv$$
$$+ \int_a^c \int_a^c \left\{ n^{-1} \sum_{i=1}^{m} \sum_{l=1}^{n_i} \sum_{r=1, r \neq l}^{n_i} E[w_{il} X_{il}^T \Sigma X_{ir} w_{ir} | U_{il} = v_1, U_{ir} = v_2] \right\}$$
$$\times f(v_1) f(v_2) \, dv_1 \, dv_2.$$

It follows, from Abel's transform, that

$$I_4(u) = K\left(\frac{b-u}{h}\right) \zeta_n - \sum_{r=1}^{n-1} \left[ K\left(\frac{d_{r+1}-u}{h}\right) - K\left(\frac{d_r-u}{h}\right) \right] \zeta_r$$

and

$$\left\| \sum_{r=1}^{n-1} \left[ K\left(\frac{d_{r+1}-u}{h}\right) - K\left(\frac{d_r-u}{h}\right) \right] [\zeta_r - n^{1/2} W(G(d_r))] \right\|_\infty$$
$$\leq \left\| \max_{1 \leq r \leq n} |\zeta_r - n^{1/2} W(G(d_r))| \sum_{r=1}^{n-1} \left| K\left(\frac{d_{r+1}-u}{h}\right) - K\left(\frac{d_r-u}{h}\right) \right| \right\|_\infty$$
$$= O_P(n^{1/4} \log n).$$

Hence,

$$I_4(u) = n^{1/2} K\left(\frac{b-u}{h}\right) W(G(b))$$
(A.12) $$- n^{1/2} \sum_{r=1}^{n-1} \left[ K\left(\frac{d_{r+1}-u}{h}\right) - K\left(\frac{d_r-u}{h}\right) \right] W(G(d_r))$$
$$+ O_P(n^{1/4} \log n)$$

uniformly for $u \in [a, b]$.

For a Wiener process, it is known that (Csörgö and Révész [5], page 44)

$$\sup_{t \in [a,b]} |W(G(t+\delta)) - W(G(t))| = O(\{\delta \log(1/\delta)\}^{1/2}) \quad \text{a.s.},$$

when $\delta$ is any small number. Using this property and the boundness of $K(\cdot)$, we have

$$\sum_{r=1}^{n-1}\left[K\left(\frac{d_{r+1}-u}{h}\right)-K\left(\frac{d_r-u}{h}\right)\right]W(G(d_r))$$

$$=\int_a^b W(G(v))\,dK\left(\frac{v-u}{h}\right)+O_P\left(\{n^{-1}\log n\}^{1/2}\right)$$

uniformly for $u\in[a,b]$. Together with (A.11) and (A.12), it follows that

(A.13)
$$\left\|(nh)^{-1/2}I_3(u)-h^{-1/2}\int_a^b K\left(\frac{v-u}{h}\right)dW(G(v))\right\|_\infty$$
$$=O_P((nh^2)^{-1/4}\log n+(nh^3)^{-1/2}).$$

Let

$$Y_{1n}(u)=h^{-1/2}\int_a^b K\left(\frac{v-u}{h}\right)dW(G(v)),$$

$$Y_{2n}(u)=h^{-1/2}\int_a^b K\left(\frac{v-u}{h}\right)[\tau(v)f(v)]^{1/2}\,dW(v-a)$$

and

$$Y_{3n}(u)=h^{-1/2}\int_a^b K\left(\frac{v-u}{h}\right)dW(v-a).$$

For a Gaussian process, following the similar proof of lemmas in Härdle [13], we have

(A.14)
$$\|Y_{1n}(u)-Y_{2n}(u)\|_\infty=O_P(h^{1/2}),$$
$$\|(f(u)\tau(u))^{-1/2}Y_{2n}(u)-Y_{3n}(u)\|_\infty=O_P(h^{1/2}).$$

Therefore, by (A.13) and (A.14),

$$\|\{nhf(u)\tau(u)\}^{-1/2}I_3(u)-Y_{3n}(u)\|_\infty=O_P((nh^2)^{-1/4}\log n+(nh^3)^{-1/2}+h^{1/2}).$$

From Theorem 2 and Theorem 3.1 of Bickel and Rosenblatt [1], when $h=n^{-\rho}, 1/5\leq\rho<1/3$, we have

$$P((-2\log\{h/(b-a)\})^{1/2}\{\nu_0^{-1/2}\|\{nhf(u)\tau(u)\}^{-1/2}I_3(u)\|_\infty-\omega_n\}<x)$$
$$\to\exp\{-2\exp\{-x\}\},$$

where $\omega_n$ is defined in Theorem 6, as

(A.15) $\qquad\operatorname{var}\{\hat\alpha_{kj}(u)|\mathcal{D}\}=\{nhf(u)\tau(u)^{-1}\nu_0^{-1}\}^{-1}(1+o_P(1))$

uniformly for $u \in [a, b]$ by straightforward calculation and similar arguments as Lemma 1. Using the same proof of the first part of Theorem 2 of Fan and Zhang [10], the result of Theorem 6 is easily obtained. $\square$

PROOF OF THEOREM 7. To prove the theorem, we first derive the rate of convergence for the bias and variance estimators. By (A.7) and its similar arguments, we have

$$\|\widehat{\text{bias}}(\hat{\alpha}_{kj}(u)|\mathcal{D}) - \text{bias}(\hat{\alpha}_{kj}(u)|\mathcal{D})\|_\infty = O_P(h^2\{n^{-1/7} + o(h)\}),$$

where the rate $n^{-1/7}$ comes from the pilot estimation of $\hat{\tilde{\boldsymbol{\theta}}}(\cdot)$ and the term $o(h)$ comes from the coefficient in front of $\hat{\boldsymbol{\theta}}^{(3)}(\cdot)$.

Furthermore, by similar proof to Lemma 1, we get

$$\left\|\frac{h}{n}H^{-1}X^TW^2XH^{-1} - f(u)\tilde{S}(u)\right\|_\infty = o_P(1)$$

where $\tilde{S}(u) = \begin{pmatrix} \nu_0 & 0 \\ 0 & \nu_2 \end{pmatrix} \otimes \Omega(u)$ and

$$\left\|\frac{h}{n}H^{-1}X^TW\mathbf{xx}^TWXH^{-1}\right\|_\infty = O_P(1).$$

These results, together with (A.7) and the results of Theorem 3 and 4, give us

$$\|\widehat{\text{var}}(\hat{\alpha}_{kj}(u)|\mathcal{D}) - \text{var}(\hat{\alpha}_{kj}(u)|\mathcal{D})\|_\infty = O_P((nh)^{-1}\{n^{-1/2} + (nh^8)^{-1/2}\}).$$

Using (A.15) and the same proof of the second part of Theorem 2 of Fan and Zhang [10], the result of Theorem 7 is obtained. $\square$

WENYANG ZHANG
DEPARTMENT OF MATHEMATICAL SCIENCES
UNIVERSITY OF BATH
UNITED KINGDOM
E-MAIL: W.Zhang@bath.ac.uk

JIANQING FAN
DEPARTMENT OF OPERATION RESEARCH
   AND FINANCIAL ENGINEERING
PRINCETON UNIVERSITY
USA
E-MAIL: jqfan@princeton.edu

YAN SUN
SCHOOL OF ECONOMICS
SHANGHAI UNIVERSITY OF FINANCE
   AND ECONOMICS
P.R. CHINA
E-MAIL: yansun2002cn@yahoo.com.cn